\documentclass{article}
\usepackage{graphicx} 

\title{CR_analysis}
\author{carolin.mehlmann}
\date{12.11.25}

\usepackage{lscape,textcomp}
\usepackage[T1]{fontenc}
\usepackage{xparse}
\usepackage{amsmath,amssymb,dsfont,units,booktabs}
\usepackage{amsthm}
\usepackage{amsfonts,import,mathtools}

\usepackage{graphicx, color}
\usepackage{color,hyperref}
\usepackage[utf8]{inputenc} 
\usepackage{xcolor}

\newtheorem{definition}{Definition}[section]

\newtheorem{remark}[definition]{Remark}
\newtheorem{lemma}[definition]{Lemma}

\newtheorem{theorem}[definition]{Theorem}

\newcommand{\nt}{\mathbf{n}}
\newcommand{\taut}{\boldsymbol{\tau}}

\newcommand{\Gr}{\color{black}}
\newcommand{\vt}{\mathbf{v}}
\newcommand{\xt}{\mathbf{x}}
\newcommand{\ft}{\mathbf{f}}
\newcommand{\gt}{\mathbf{g}}
\newcommand{\ut}{\mathbf{u}}
\newcommand{\wt}{\mathbf{w}}
\newcommand{\zt}{\mathbf{z}}

\newcommand{\Vt}{\mathbf{V}}

\newcommand{\Pt}{\mathbf{P}}
\newcommand{\It}{\mathbf{I}}

\newcommand{\phit}{\pmb{\phi}}
\renewcommand{\div}{\operatorname{div}}
\newcommand{\tr}{\operatorname{tr}}
\newcommand{\et}{\mathbf{e}}

\begin{document}

\title{Analysis of the Crouzeix-Raviart Surface Finite Element Method for Vector-Valued Laplacians}

\author{\underline{Carolin Mehlmann}\thanks{Institute of Analysis and Numerics,
    Otto-von-Guericke University Magdeburg, Universitätsplatz 2,
  39106 Magdeburg, Germany, carolin.mehlmann@ovgu.de}}

\maketitle

\begin{abstract}
  Recently, a nonconforming surface finite element was developed to discretize 3d vector-valued compressible flow problems arising in climate modeling. In this contribution we derive an error analysis for this approach on a vector-valued Laplace problem, which is an important operator for fluid-equations on the surface. In our setup, the problem is approximated via edge-integration on local flat triangles using the nonconforming linear Crouzeix-Raviart element. The latter
   is continuous at the midpoints of the edges in each vector component. 
This setup is numerically efficient and straightforward to implement. For this Crouzeix-Raviart discretization we introduce interpolation estimates, derive optimal error bounds in the $H^1$-norm and $L^2$-norm and present an estimate for the geometric error. Numerical experiments validate the theoretical results. 
\end{abstract}

\section{Introduction}
\subsection{Model problem and earlier work}
Low order nonconforming finite element discretizations play an important role for the approximation of geophysical flow problems on the surface of the Earth \cite{MehlmannGutjahr2022, COMBLEN2009}. On the used meshes in climate models this type of discretization is  a good compromise between accuracy and efficiency of the numerical setup \cite{Mehlmann2021,Mehlmannetal2021,COMBLEN2009}. In this paper, we analyze a recently developed nonconforming Crouzeix-Raviart approximation \cite{MehlmannGutjahr2022}, which has also been studied in \cite{Mehlmann2023Pamm}. The nonconforming discretization has been developed for the approximation of sea-ice dynamics in the climate model ICON \cite{JUNGCLAUS_2020}.

An important aspect in the analysis of finite element methods for flow problems on surfaces is to ensure that the numerical approximation of the flow is tangential to the surface. This difficulty already arises in the analysis of the surface Laplace equation.
Vector-valued Laplace problems on the surface are easier to analyze than the Stokes problem because the equation consists only of the velocity vector and no pressure unknown.
Therefore, vector-valued Laplace equations are a useful simplification of more complex surface flow problems for finite element analysis, which is the topic treated in this paper. 

We consider an oriented, connected, bounded, $C^\infty$-smooth surface $\Gamma$ in $\mathbb{R}^3$ with $\partial \Gamma=\emptyset$. The vector-valued problem based on the Bochner Laplace is given as 
\begin{equation}\label{eq:conEq}
    -\Delta_\Gamma \ut +\ut =\ft.
\end{equation}
The zero order term $\ut$ has been added to the left-hand side of the equation to avoid technical details  related to the kernel of the Laplace operator.
 
In the work of \cite{Gross2018,Jahnkuhn2020,HansboLarson2020, Hardering2022}  finite element methods are 
studied for a surface vector Laplace problem and  error bounds are derived. The authors consider $H^1$-conforming finite elements and enforce tangentiality condition of the vector-valued velocity field either by a penalty term or a Lagrange multiplier. Other approaches to the surface Stokes problem abandon the $H^1$-conformity but enforce tangentiality through the construction of the elements (e.g., \cite{Bonito2020e,demlow2023tangential}). The Crouzeix-Raviart approach presented here is likewise $H^1$-nonconforming and, by construction, tangential to the discrete surface. The element can be placed within the framework of surface finite elements (SFEM), where the surface is approximated by a geometry $\Gamma_h$, usually a polyhedron. {We establish in each face of the polyhedron a local basis and express the Crouzeix-Raviart approximation in terms of the conormal vectors $\nt_E,\taut_E$ given in the edge midpoint. While $\taut_E$ is continuous, $\nt_E$ is discontinuous along an edge.  A sketch of the setting is shown in Figure \ref{fig:set_up}. The discretization is straightforward to implement and benefits from its efficiency; see for instance \cite{Mehlmannetal2021}.

To our knowledge, this paper is the first to present a rigorous mathematical formulation and analysis of a vector-valued surface Crouzeix-Raviart element, including optimal error estimates in both the $L^2$- and $H^1$-norms. An analysis of the scalar surface Crouzeix-Raviart element has previously been given in \cite{Guo2020}.
Compared to existing approaches, the proposed Crouzeix-Raviart discretization performs better in the linear setting than existing conforming methods (e.g., \cite{HansboLarson2020}, \cite{Hardering2022}) without additional geometric information and  other penalty-free approaches evaluated on the surface Laplace, such as DG methods in \cite{Lederer2019}. A similar penalty-free approach has been derived for the incompressible Stokes equation (e.g. \cite{demlow2023tangential}), but the work does not present an $L^2$-error estimate.



\subsection{Main contributions}
This section provides an informal summary of the paper’s main contributions. Technical details, including the formal definition of operations on functions defined over different domains, are postponed to Section~\ref{sec:2}.

\paragraph{A-priori energy estimate}
With respect to the exact geometry \(\Gamma\), we define the  bilinear form \(a(\cdot,\cdot):=(\nabla_\Gamma\, \cdot,\nabla_\Gamma\, \cdot) + (\cdot,\cdot)\)  to write~\eqref{eq:conEq} in variational formulation  
\[
a(\ut,\vt) = (\ft,\vt)_\Gamma\quad \forall \vt \in H^1_{tan}(\Gamma),
\]  
where \(\ut\in H^1_{tan}(\Gamma)\) denotes the exact solution and $H^1_{tan}(\Gamma)$ is the Sobolev space of the tangential vector fields with covariant derivatives. { We extend functions $\vt$ from $\Gamma$ to $\Gamma_h$ along the normal directions. The extension is denoted by $\tilde \vt$.  }
The corresponding discrete weak formulation reads as follows.
\begin{align}
a_h(\ut_h,\vt_h)=(\tilde \ft,\vt_h)_{\Gamma_h} \quad \forall \vt_h \in \Vt_h,
\end{align}
where $\Vt_h$ is the space of the vector-valued surface Crouzeix-Raviart functions on $\Gamma_h$ and $\ut_h$ is the finite element solution. 
Here, \(a_h(\cdot,\cdot)\) denotes the symmetric bilinear form of the Crouzeix-Raviart approximation integrated on $\Gamma_h$, formulated using approximate surface differential operators based on \(\Gamma_h\). For this setting, the following error estimate (see Theorem~\ref{th:energy}) holds.
 \begin{equation}\label{eq_energy1}
  \|\tilde \ut-\ut_h\|_h \leq c h \|\ft\|_{L^2(\Gamma)},
  \end{equation}
  where $\| \cdot\|_h$ is the energy norm and $c$ a mesh independent constant. The proof of this estimate is based on the second Strang Lemma:
 \begin{equation}
  \|\tilde \ut-\ut_h\|_h \leq \underbrace{c \inf_{\vt_h \in \Vt_h} \|\tilde \ut-\vt_h\|_h}_{\text{approximation error}} + \underbrace{\sup_{\vt_h \in \Vt_h} \frac{a_h(\tilde \ut, \vt_h) - (\tilde \ft,\vt_h)}{\| \vt_h\|_h}}_{\text{nonconformity error}}.
\end{equation}

To bound the \emph{approximation error}, we establish a set of interpolation estimates in Lemma~\ref{lemma_INT} for the interpolant  $\Pi^{{tan}}_h : H^1_{{tan}}(\Gamma) \to {\Vt}_h$. The key idea of the proof is to introduce a vector-valued extension of the scalar Crouzeix-Raviart interpolant, $\Pi_h\vt \not\subset \Vt_h$, and to derive estimates for the difference between the two interpolations, $\Pi^{{tan}}_h$ and $\Pi_h$ . Due to the nonconformity of the Crouzeix-Raviart element, we estimate the \emph{nonconformity error} by working with the strong formulation and applying integration by parts. This gives rise to the typical jump terms across element edges. The main challenge in estimating these cross-element jump terms lies in controlling the components of the conormal vectors, since the corresponding conormals \(\mathbf{n}_E\) of two triangles sharing an edge do not lie in a  plane (see Figure~\ref{fig:set_up}, panel (c)). This is a distinctive feature of the vector-valued surface Crouzeix-Raviart finite element (\ref{eq:VT}). A bound for the nonconformity error is shown in Lemma \ref{lemma:noncon}. The geometric errors arising in the nonconformity estimate can be bounded by the arguments provided in \cite{HansboLarson2020}.
    
\paragraph{Analysis for the $L^2$-estimate}
To obtain an $L^2$-estimate, we consider the projected error on the discrete surface $\Pt_h(\tilde \ut- \ut_h)$ in the $L^2$-norm. This is due to the fact that $\tilde \ut -\ut_h=\Pt \tilde \ut-\Pt_h\ut_h$ is only $O(h)$ as $\ut$ and $\ut_h$ lie in different tangential planes. Since Galerkin orthogonality fails in the presence of the geometric approximation due to $\vt_h \in \Vt_h\not \subset H^1(\Gamma)$ we have $a_h(\tilde \ut-\ut_h,\vt_h) \not= 0$, and an additional remainder terms arise when applying the a duality argument to obtain an $L^2$-estimate (see Theorem \ref{th:L2}). In the proof of the $L^2$-estimate, we obtain the expression 
\begin{align*}
  \| \Pt_h(\tilde \ut-\ut_h)\|^2_{L^2(\Gamma_h)}
        = & [(\tilde \ut,\tilde \gt)_{\Gamma_h} -(\ut,\gt)_\Gamma]+ [a_h(\tilde \ut,\tilde \zt)-a(\ut,\zt)]\\
        &+  a_h(\tilde \ut-\ut_h, \tilde \zt-\zt_h)\\
       & +[(\ut,\gt)_\Gamma-(\ut_h,\tilde \gt )-a_h(\tilde \ut-\ut_h,\tilde \zt)]\\
        &+[(\ft,\zt)_\Gamma-(\tilde \ft,\zt_h)-a_h(\tilde \ut,\tilde \zt -\zt_h)],
    \end{align*} 
    where the first  and second term is {$O(h^2)$} due to geometric error estimates provided in \cite{HansboLarson2020}, the third term is also $O(h^2)$ due to the interpolation estimates we derived  in the energy norm. The fourth and fifth terms are also of order $O(h^2)$ according to the primal and dual consistency estimate provided in Lemma \ref{lemmma:dual_consist}. The latter estimate is proven by inserting and subtracting the interpolant $\Pi^{tan}_h$ and by deriving the primal and dual interpolation estimates in Lemma \ref{lemma:Klaus2}.

\subsection{Outline}

The remainder of paper is structured as follows. Section~\ref{sec:2} begins with the geometric description of the surface $\Gamma$, followed by the definition of a mapping between functions on $\Gamma$ and $\Gamma_h$. We then derive the corresponding surface differential operators, which provide the foundation for introducing Sobolev spaces of arbitrary order on curved surfaces. In Section~\ref{sec:geo}, we introduce the surface approximation and present the geometry results for a triangulation $\Gamma_h$ of the exact surface $\Gamma$. In addition, we establish norm equivalences for surface Sobolev norms between functions on $\Gamma$ and on $\Gamma_h$. In Section~\ref{sec:nonconform}, we first introduce the weak formulation of the Bochner–Laplace problem on surfaces.
We then present the surface vector-valued Crouzeix-Raviart method and formulate the discrete weak problem using surface differential operators and integration on the approximate surface $\Gamma_h$.  The error analysis is carried out in Section \ref{sec:error}. The main results of the section are interpolation estimates, and the optimal error bounds in $H^1$-norm and $L^2$-norm. Numerical examples to validate the theoretical results are presented in Section \ref{sec:num}.

\section{Surface, derivatives and norms}\label{sec:2}
The signed distance function of $\Gamma$ is given by $d(x)$.  
Using the distance function, we define the outward pointing unit normal vector  as \begin{equation}
\nt(x)=\nabla d(x),
\end{equation}
where $\nabla$ is the standard gradient in  $\mathbb{R}^3$.
We introduce a neighborhood, a strip around $\Gamma$ with distance $\delta$, as
\begin{equation}\label{eq:U}
U=\{ x \in\mathbb{R}^3 : dist(x,\Gamma)  \leq \delta \},
\end{equation}
where $dist(x,\Gamma)$ is the Euclidean distance between $x$ and $\Gamma$. Let $\delta$ be small enough such that a unique closest point mapping $p(x)$ from $U \to \Gamma$ exists  with
\begin{equation}\label{eq:proh}
    p(x)=x-d(x)\nt(p(x)).
\end{equation}
Let $\phi: \Gamma \to \mathbb{R}$. By $\tilde \phi:U\to \mathbb{R}$ we denote the extension to the neighborhood $U$, defined along the normal directions as 
\begin{equation}\label{eq:extension}
    \tilde \phi(x)=\phi(p(x)), \, \forall x \in U.
\end{equation}
{Later on $\Gamma_h$ we will denote the discrete surface which is defined by the triangulation and $\Gamma$ will represent the lifted counterpart. The usage will be clear from the given context.} Analogously, the lift of  a function $\phi:\Gamma_h \to \mathbb{R}$ on $\Gamma_h \subset U$ to $\Gamma$ is given as  \begin{equation}
    \phi^l(x)=\phi(\eta(x)), \, \forall x \in \Gamma,
\end{equation}
where $\eta(x)$ is the unique solution of 
\begin{equation}
   x=p(\eta)=\eta-d(\eta)n(x).  
\end{equation} 
Throughout the paper, we apply a component-wise lifting and extension of vector-valued functions.
\subsection{Tangential differential operators}
To define the surface derivatives, we introduce the projection into the tangential plane as
\begin{equation}\label{eq:P}
 \Pt=\It-\nt 
   \nt^T.
   \end{equation}
 The covariant derivatives of a vector field {$\vt$} are defined as 
\begin{equation}
  \nabla_\Gamma \vt:= \Pt \nabla \tilde \vt \Pt.
\end{equation}
Similarly to \cite{Gross2018} we  introduce the divergence operator  for vector-valued ${  \vt}: \Gamma \to \mathbb{R}^3$ and matrix valued objects $A:\Gamma \to  \mathbb{R}^3 \times  \mathbb{R}^3$ as 
{
\begin{align*}
    \div_\Gamma \vt&=\tr(\nabla_\Gamma  \vt)=\tr(\Pt \nabla \tilde \vt \Pt),\\
     \div(A):&=\Big( \div_\Gamma(\et_1 ^ T A), \div_\Gamma(\et_2^ T A), \div_\Gamma(\et_3^ T A) \Big ),
\end{align*}
where $\et_1,\et_2,\et_3$ are the unit vectors of $\mathbb{R}^3$.}
The surface Laplacian is defined by
\begin{align*}
 \Delta_\Gamma \vt:&=\Pt\div_\Gamma(\nabla_\Gamma  \vt).
\end{align*}
 This is the so called Bochner Laplace, which is also treated in \cite{HansboLarson2020}. In differential geometry and in exterior calculus, a different Laplace operator, the Hodge Laplace, is considered. Another surface Laplace operator based on the symmetric gradient $\nabla \vt +\nabla \vt^T$ is analyzed in \cite{Gross2018}.
  
\subsection{Surface Sobolev space}
{Let $\omega \subset \Gamma$. In the following, we introduce the surface Sobolev space of $k$ times  weakly differentiable vector valued functions with  componentwise derivatives } and the corresponding norm as 
{
\begin{equation}\label{eq:sobh}
    \|\vt\|^2_{H^k(\omega)} := \sum^k_{j=0}\|(\nabla \Pt)^j \tilde \vt\|^2_{L^2(\omega)},\, |\vt|^2_{H^k(\omega)} :=\|(\nabla \Pt)^k \tilde \vt\|^2_{L^2(\omega)}.
 \end{equation}
 } 
 We write  $H^k(\Gamma):=\Big (H^k(\Gamma)\Big )^3$ for vector-valued functions  in $\mathbb{R}^3$.
 The space of the vector-valued functions that are tangential to the surface is named as 
\begin{equation}
    H_{tan}^k(\Gamma):=\{\vt \in H^k(\Gamma)|\vt \cdot \nt=0\}.
\end{equation}

The corresponding norm is given as
\begin{equation}\label{eq:sobh_tan}
   \|\vt\|^2_{H^k_{tan}(\omega)} := \sum^k_{j=0}\|(\nabla_\Gamma)^j \vt\|^2_{L^2(\omega)}.
\end{equation}
As shown in Lemma 2.3 of \cite{HansboLarson2020}, for fields $\vt \in H^k_{\mathrm{tan}}(\Gamma)$, the norms $\|\vt\|_{H^k(\Gamma)}$ and $\|\vt\|_{H^k_{\mathrm{tan}}(\Gamma)}$ are equivalent.

\section{Geometry approximation}\label{sec:geo}
\subsection{Approximation of the surface}
We denote by $\Gamma_h=\cup_{i\in \mathbb{N}} K_i\subset U$ the triangulation of $\Gamma$. The triangulation is shape regular and quasi uniform with the maximal diameter 
$
h=\max_{K\in \Gamma_h} \operatorname{diam}(K).
$
The triangulation is made such that all {vertices $\xt$} lie on $\Gamma$.  
\begin{figure}
\begin{center}
    \begin{tabular}{l l l l}
    a)& &b)&c)\\
   \includegraphics[width=0.29\textwidth]{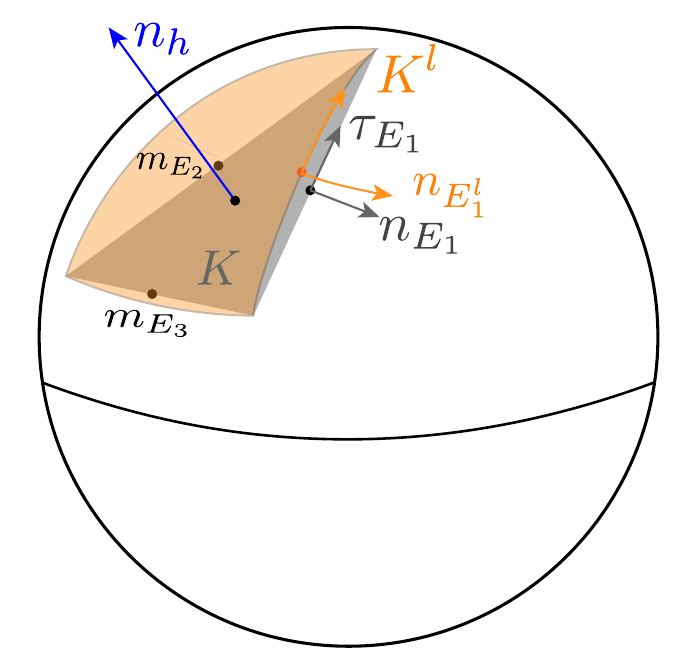}  &  &\includegraphics[width=0.29\textwidth] {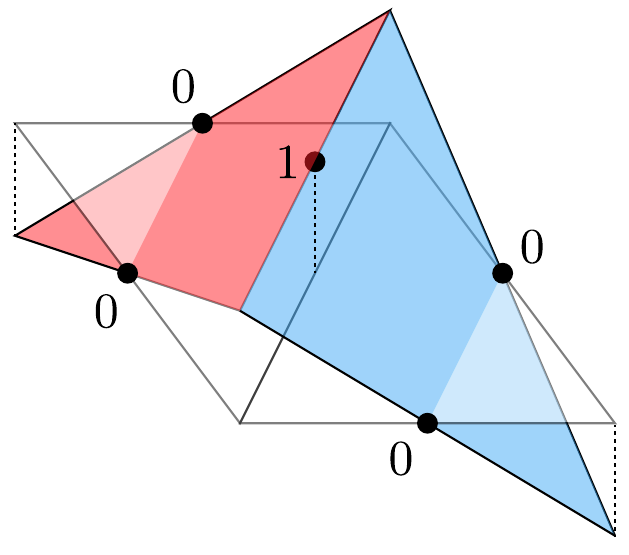} &
   \includegraphics[width=0.29\textwidth]{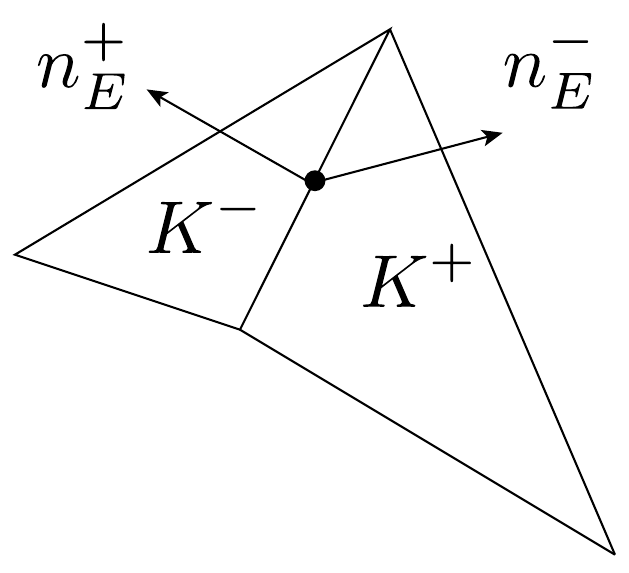}
\end{tabular}
\end{center}
\caption{{a):} The tangential plane which coincides with the sphere at the vertices of the triangle. The outward pointing normal vector (blue) is named {$\nt_h$.} 
The tangential and conormal of the flat triangle $K$ (shaded in gray) to an edge $E$ is called {$\nt_E, \taut_E$}. 
 Accordingly the tangential and conormal of the curved  triangle $K^l$ (orange) is named $\nt_E^l, \taut_E^l$. {b)}: Crouzeix-Raviart basis function. {c)}: The outward pointing conormal vectors $\nt_E^\pm$ that share an edge $E$. 
 \label{fig:set_up}}
\end{figure}
Each triangle $K$ consists of three vertices $\xt_{1}, \xt_{2},\xt_{3}$ numbered counterclockwise  and edges $E_{1}=(\xt_1,\xt_2), E_2=(\xt_2,\xt_3), E_3=(\xt_3,\xt_1)$. Based on the vectors $\vec{E_1}=\xt_2-\xt_1$, $\vec{E_2}=\xt_3-\xt_2$, $\vec{E_3}=\xt_1-\xt_3$ the outward pointing normal vector is calculated by
\begin{align*}
 \nt_h^K= \frac{\vec{E_1}\times (-\vec{E_3})}{\|\vec{E_1}\times (-\vec{E_3})\|}.
\end{align*} 
Using $\nt^K_h$ the projection onto the tangent space of $\Gamma_h$ is defined by
\begin{equation}
    \Pt_h=\It-\nt^K_h (\nt^K_h)^T.
\end{equation}
The setup is visualized in the left panel of Figure \ref{fig:set_up}. 
The conormal vectors to an edge $E_i, \,i=1,..,3$ are given by
\begin{align*}
\taut^K_{E_i}=\frac{\vec{E}_{i}}{\|\vec{E}_{i}\|}, \quad \nt^K_{E_i}=\taut^K_{E_i }\times \nt^K_h, \quad i=1,2,3.
\end{align*}  
To any  $K \in \Gamma_h$ a curved triangle $K^l=p(K)$ on $\Gamma$ is related such that  the set of curved triangular faces is defined as $\Gamma^l_h=\{K^l: K \in \Gamma_h \}$ with $
    \Gamma= \cup_{K^l \in \Gamma^l_h} K^l.$ Any edge midpoint $m_E$ of edge $E$ is shared by two triangles $K^{+}$ and $K^{-}$.  The set of all edge midpoints is denoted by ${\cal{M}}_h$. We denote the outward pointing conormal vectors as $\nt_E^\pm$. Note that $\nt^+_E\not= \nt_E^-$, but $\taut^+_E=-\taut^-_E$. The conormal vectors of a curved edge $E^l=p(E)$ are named $\nt^{\pm}_{E^l}, \taut^\pm_{E^l}$. A sketch of the setup is shown in the right panel of Figure \ref{fig:set_up}. 
    
\subsection{Geometric approximation results}    
{ We introduce the geometric error bounds. For the sake of readability, we will omit the tilde when indicating the extensions of $\Pt$ and $\nt_{E^l}$ and the subscript $l$ for the lift, whenever this is clear from the context.} The following geometric approximation holds true for our setup.
\begin{lemma}[Geometry approximation]
Let $\Gamma_h \subset U$ be an approximation  of $\Gamma$ with properties outlined above. Assume that the mesh size $h$ is small enough. Then it holds that 
\begin{align}
           \|\Pt-\Pt_h\|_{L^\infty(K)}&\leq c h\label{eq:P_Ph},\\
           \|\nt^\pm_{E^l}-\nt_E^\pm\|_{L^\infty(E)}&\leq c h, \label{eq:nel}\\
           \|\nt^{\pm}_{E^l}-\Pt \nt^\pm_E\|_{L^\infty(E)}&\leq c h^2,\label{eq:pm_n}\\
            \|\Pt_h\nt^{\pm}_{E^l}-\nt^\pm_E\|_{L^\infty(E)}&\leq c h^2.\label{eq:pm_n2}
\end{align}
\end{lemma}
\begin{proof}
 The proof of (\ref{eq:P_Ph})-(\ref{eq:pm_n}) is given in Appendix A of \cite{Larsson2013}. The estimate \eqref{eq:pm_n2}
is shown by noting that 
\begin{multline*}
   \|\Pt_h\nt^{\pm}_{E^l}-\nt^\pm_E\|_{L^\infty(E)} = \|\Pt_h(\nt^{\pm}_{E^l}-\nt^\pm_E)\|_{L^\infty(E)}\\
  \leq  \|\Pt (\nt^{\pm}_{E^l}-\nt^\pm_E)\|_{L^\infty(E)}+\|(\Pt_h - \Pt ) (\nt^{\pm}_{E^l}-\nt^\pm_E)\|_{L^\infty(E)}
\end{multline*}
and  applying \eqref{eq:P_Ph}-\eqref{eq:pm_n}. 
\end{proof}    

\subsection{Sobolev norm approximation}
Since the discrete surface $\Gamma_h$ is only piecewise smooth, we introduce the {broken Sobolev space} $H^k(\Omega_h)$ and $H^k_{tan}(\Omega_h)$ on $\Gamma_h$, defined for  vector-valued functions. The associated norms are

\begin{equation}
    \|\vt\|^2_{H^k(\Omega_h)} := \sum^k_{j=0}\|(\nabla \Pt_h)^j \vt\|^2_{L^2(\Omega_h)}, \quad \|\vt\|_{H_{tan}^k(\Omega_h)}^2 
  = \sum_{j=0}^k \| (\nabla_{\Gamma_h})^j \vt \|_{L^2(\Omega_h)}^2.
  \label{eq:broken_sobolev_norm}
 \end{equation}
These definitions are analogous to (\ref{eq:sobh}) and (\ref{eq:sobh_tan}), but applied to the discrete surface $\Gamma_h$. 
Whenever $\Omega_h$ is used as the integration domain, the integration is carried out 
elementwise,
\[
  \| \cdot \|_{L^2({\Omega_h})}^2 = \sum_{K \in \Omega_h} \| \cdot \|_{L^2(K)}.
\]

The norms on $\Gamma$ and $\Gamma_h$ are related by the following bounds:
\begin{lemma}[Norm equivalence]\label{lemma:norm}
Let $K\in \Gamma_h$ and $K^l\in\Gamma$ be its lifting.
 If $\vt \in H^2(K^l)$, then for any $K^l\in \Gamma$ it holds that
\begin{align}
   \|\vt\|_{L^2(K^l)}&\leq c \|\tilde \vt\|_{L^2(K)} \leq c\|\vt\|_{L^2(K^l)},\label{eq:norm_l2}\\
   |\vt|_{H^1(K^l)} &\leq c |\tilde \vt|_{H^1(K)} \leq c |\vt|_{H^1(K^l)},\label{eq:norm_equ}\\
    |\tilde \vt|_{H^2(K)} &\leq c |\vt|_{H^2(K^l)},\label{eq:norm2}\\
   |\vt|_{H^2(K^l)}& \leq c |\tilde \vt|_{H^2(K)} \label{eq:norm3}.
\end{align}
\end{lemma}
\begin{proof}
Inequalities (\ref{eq:norm_l2})-(\ref{eq:norm_equ}) are proven in \cite{HansboLarson2020} in Section 4.2.
Inequalities (\ref{eq:norm2}) and  (\ref{eq:norm3}) are proven for the scalar case in Appendix B of \cite{Larsson2013} and can be directly translated to the vector-valued setup.
\end{proof}

\section{The nonconforming finite element approximation on the surface}\label{sec:nonconform}
In this section, we consider the Bochner-Laplace problem on the surface $\Gamma$ and derive its weak formulation. We then introduce the vector-valued surface Crouzeix-Raviart problem, formulated on the approximated surface $\Gamma_h$ using the corresponding surface differential operators based on $\Gamma_h$.

\subsection{Weak form of the Bochner-Laplace problem}
The corresponding weak formulation of the Bochner-Laplace problem (\ref{eq:conEq}) reads as: 
\begin{equation}\label{eq:ana}
 \text{Find} \quad \ut \in H_{tan}^1(\Gamma)\, \text{s.t.} \quad   a(\ut,\vt)=(\ft,\vt)\, \,\forall \vt \in H_{tan}^1(\Gamma), 
\end{equation}
with $a(\vt, \ut)=  (\nabla_\Gamma \ut, \nabla_\Gamma \vt  )_\Gamma + ( \ut, \vt )_\Gamma, \, \ut, \vt \in  H_{tan}^1(\Gamma)$, where $(\cdot, \cdot)$ is the standard $L^2$-inner product, and  $f \in L^2(\Gamma)$. The Bochner-Laplace equation is analyzed in \cite{HansboLarson2020}. In the study, the authors demonstrated by the Lax-Milgram theorem the uniqueness of the solution of equation (\ref{eq:ana}) and derive the following regularity estimates
\begin{equation}\label{eq:regu}
    \|\ut\|_{H^2_{tan}(\Gamma)} \leq c \|\ft\|_{L^2(\Gamma)},
\end{equation}
where $c$ is a generic positive constant.

\subsection{The Crouzeix-Raviart Finite Element space}
In the following, we construct the finite element space of the vector-valued Crouzeix-Raviart functions that are tangential to $\Gamma_h$.
We consider a uniform reference triangle  $\hat K$ with edges of length $1$. $\hat K$ is related to each $K \in \Gamma_h$ via an affine linear mapping $F_K: \hat K \to K$.
{In each edge midpoint $m_{\hat E_i}$, with  $i=1,2,3$, of the reference triangle $\hat K$ the linear Crouzeix-Raviart basis functions are given by the defining properties 
\begin{equation}\label{phi}
\hat \phi_i \in P^1(\hat K), \quad \hat \phi_i(m_{\hat E_j})=\delta_{ij}, \quad  \forall i,j=1,\dots,3.
\end{equation}
Hence, $\int_{\hat E_j} \hat \phi_i d\hat s=\delta_{ij}, \,\forall i,j=1, \dots,3.$
The coupling of two Crouzeix-Raviart basis functions on two neighboring elements is shown in the middle panel in Figure \ref{fig:set_up}.

We indicate a finite element function locally by its conormal components in the planes spanned by each element $K$:
\begin{equation}
\Vt_K:= \text{span}\langle \nt^K_{E_{i}}\phi^K_{i},\taut^K_{E_{j}}\phi^K_{j},\; i,j=1,2,3\rangle,
\end{equation}
where $\phi^K_{i}=\hat \phi_i \circ F_K^{-1}$. Thus, we get 
 \begin{align*}
 \vt_h|_K&=\sum^3_{i=1}(v^{n,K}_{i} \nt^K_{E_i}+v^{t,K}_{i}\taut^K_{E_i})\phi^K_{i}
\end{align*}
with coefficients $v_i^{n,K}\in\mathbb{R}$ and $v_i^{t,K}\in\mathbb{R}$ for  $i=1,2,3$. 
Based on \eqref{phi} it holds that 
 \begin{align}
{ \partial_{\tau_{E_i}} \phi^K_i|_{E_i}= 0,\ \quad  \forall i=1, \dots,3.}
\end{align}  
The global space of these vectors  is given by
\begin{align*}
\bar \Vt_h=\Big \{ \vt: \Gamma_h \to \mathbb{R}^3, \vt|_K \in \Vt_K \}.
\end{align*}
By construction, all $\vt_h\in \bar \Vt_h$ satisfy $\vt_h \cdot \nt_h=0$.
Using this construction, we introduce the space of the tangential Crouzeix-Raviart functions on $\Gamma_h$, which are continuous for each conormal-component at the edges midpoints by
\begin{multline}\label{eq:VT}
\Vt_h=\Big\{ \vt \in \bar \Vt_h: \vt|_K \in \Vt_K,
\quad\forall E=\partial K^+\cap\partial K^-\;(K^+\neq K^-):\\
\vt\big|_{K^+}(m_E)\cdot \nt^+_{E}=-
\vt\big|_{K^-}(m_E)\cdot \nt^-_{E},\\
\vt\big|_{K^+}(m_E)\cdot \taut^+_{E}=-
\vt\big|_{K^-}(m_E)\cdot \taut^-_{E}\Big \}.
\end{multline}
Let the jump of a vector-valued function  across an edge $E$ be defined as
\begin{equation}\label{eq:jump}
    [\vt]= \lim_{s \to 0_+} \Big( \vt(x-s \nt^+_E) - \vt(x-s\nt_E^-) \Big) .
\end{equation}
Then, based on the continuity condition of the vector components at the edge midpoints  formulated in (\ref{eq:VT}) it follows that the average over the normal and  the tangential component is zero, and we get that { \begin{equation}
\int_E[\vt_h \cdot \nt_E] \,d \sigma_h=0, \quad \int_E[\vt_h \cdot \taut_E] \, d\sigma_h=0.
\end{equation}
}

{ Finally, we introduce {the space of the scalar Crouzeix-Raviart finite element, which is needed for the derivation of the interpolation estimates in Section \ref{sec:Int}. On $\Gamma_h$ it is given by 
\begin{equation}
V_h^{CR}=\{\psi \in L^2(\Gamma_h): \psi|_K \in P^1(K) \text{ and } \psi \text{ is continuous at } {\cal{M}}_h\},
\end{equation}
where ${\cal{M}}_h$ is the set of all edge midpoints.}}

 
 \subsection{Broken weak formulation}
The weak formulation of our continuous problem  (\ref{eq:ana}) reads as: 
\begin{equation}\label{eq:var}
\text{Find } \ut_h \in \Vt_{h}\text{ s.t. }  a_h(\ut_h,\vt_h)=(\tilde \ft,\vt_h)\quad \forall \vt_h \in \Vt_{h}, 
\end{equation}
with the bilinear form 
\begin{equation}\label{eq:ana_disc}
    a_h(\ut_h, \vt_h)= \sum_{K \in \Gamma_h}  ( \nabla_{\Gamma_h} \ut_h,\nabla_{\Gamma_h} \vt_h )_K + (\ut_h,\vt_h)_{\Gamma_h}.
\end{equation}
The discrete vector-valued derivative 
$\vt_h \in \Vt_{h}$ on each element $K$ is defined as 
 \begin{align}\label{eq:discret_grad}
     \nabla_{\Gamma_h}\vt_h&=\Pt_h \nabla \vt_h \Pt_h.
 \end{align} 
Based on the weak formulation, the energy norm is given by 
\begin{equation}\label{eq:enegry_norm}
    \|\vt_h \|^2_h= |\vt_h|^2_{H^1_{tan}(\Gamma_h)} + \|\vt_h\|^2_{L^2(\Gamma_h)}=a_h(\vt_h,\vt_h).
\end{equation}
We note that $\|\vt_h \|_h$ is a norm on $\Vt_h$.  
{

{The discrete variational formulation (\ref{eq:var}) is elliptic and continuous and therefore has a unique solution.} The discrete problem introduces three nonconformities compared to the original formulation (\ref{eq:ana}), namely the 
geometric  error $\Gamma_h \neq \Gamma$, the tangential inconsistency $\Pt_h \neq\Pt$ and { the $H^1$-nonconformity}.

\section{A priori error estimates}\label{sec:error}
We now establish error estimates for our method in both the energy-norm and $L^2$-norm, based on the assumptions and approximation properties presented in Section~\ref{sec:geo}. We begin by deriving the interpolation estimate in Section~\ref{sec:Int}. In Section~\ref{sec;energy}, we introduce Lemma \ref{lemma:noncon} to bound the nonconformity error caused by the use of the Crouzeix-Raviart element; the result is then employed in the proof of the main energy estimate in Theorem~\ref{th:energy}.
The $L^2$-error estimate is derived in Section~\ref{Sec:L2}, where we first present a primal-dual interpolation estimate (Lemma \ref{lemma:Klaus2}) and a primal-dual consistency error estimate (Lemma (\ref{lemmma:dual_consist})), which together form the basis of the proof of Theorem~\ref{th:L2}.

\subsection{Interpolation operator}\label{sec:Int}
{ We proceed by defining the interpolation operator $\Pi_h^{tan}: H^1_{tan}(\Gamma)\to \Vt_h$ as a Scott-Zhang interpolation operator, which is $H^1$-stable cf.~\cite[Th. (3.1)]{Scott1990}. 
On each cell $K$, let $\phi_i$, with $i=1,\dots,3,$ denote the Crouzeix-Raviart basis functions, and let $E_i,  i=1,\dots,3,$ be the edges of $K$. The vectors  $\nt_{E_i}, \taut_{E_i}$, $i=1,\dots,3,$ are the corresponding conormal vectors to the edge $E_i$. Furthermore,  $\widetilde{ \vt \cdot \nt_{E_i^l}}, \widetilde{ \vt \cdot \taut_{E_i^l}}$  denote the extension of the  two conormal components of $\vt$ from $K^l$ to $K$. Based on this notation, the interpolant is constructed as follows. }
\begin{multline}\label{eq:IntCompTan}
 \Pi_h^{tan}\vt|_K=\sum^3_{i=1}  \phi_i\alpha_i^{tan}(\tilde \vt),\\ \alpha_i^{tan}(\tilde \vt)=\frac{1}{|E_i|}\int_{E_i} (\widetilde{ \vt \cdot \nt_{E_i^l}}) \nt_{E_i} \, d\sigma_h + \frac{1}{|E_i|}\int_{E_i} (\widetilde{ \vt \cdot \taut_{E_i^l}})  \taut_{E_i} \, d\sigma_h.
\end{multline}
For $\Pi^{tan}_h$ the following properties hold:
\begin{align}\label{eq:int_tan}
\frac{1}{|E_i|}\int_{E_i}\Pi_h^{tan} \vt \cdot \nt_{E_i} \,d \sigma_h&=\frac{1}{|E_i|}\int_{E_i}\widetilde{ \vt \cdot \nt_{E_i^l}}\, d\sigma_h,\\
\frac{1}{|E_i|} \int_{E_i}\Pi^{tan}_h \vt\cdot \taut_{E_i}\, d\sigma_h&=\frac{1}{|E_i|}\int_{E_i} \widetilde{ \vt \cdot  \taut_{E_i^l}} \, d\sigma_h.
\end{align}}
This means that the components with respect to $\Gamma$ are interpolated into components with respect to $K$, see Figure \ref{fig:set_up} a) for a visualization.

{To derive interpolation estimates for $\Pi_h^{tan}$, we introduce a vector-valued extension of the scalar surface Crouzeix-Raviart interpolant \cite{Guo2020}. We start by recalling the construction of the latter:} 
Let $K \in \Gamma_h$ and $\psi$ denote a scalar function with $\psi \in H^1(K^l)$.  On $K$ the Crouzeix-Raviart interpolation $\Pi_K: { H^1(K^l) \to P^1(K) }$ is defined as 
{\begin{equation}\label{eq:interpol}
    (\Pi_K \psi)(m_{E_i})= \frac{1}{|E_i|}\int_{E_i} \tilde \psi d\sigma_h,\, i=1,2,3 \quad \forall  \psi \in H^1(K^l),
\end{equation}
where 
 $m_{E_i}$, $i=1,2,3$, are the edge midpoints of $K$.} The local interpolation estimate \cite{Crouzeix1974} is given by 
\begin{equation}\label{eq:localerror}
    \|\tilde \psi-\Pi_K \psi\|_{L^2(K)}+h| \tilde \psi-\Pi_K \psi|_{H^1(K)} \leq c h^2|\tilde \psi|_{H^2(K)},
\end{equation}
where $\tilde \psi \in H^2(K)$. 
}
 The global interpolation operator $\Pi_h:H^1(\Gamma) \to \Vt^{CR}_h$} is defined as 
\begin{align}\label{globalInt}
    (\Pi_h \psi)|_K=\Pi_K \psi, \forall K \in \Gamma_h. 
\end{align}
{For a vector-valued function $\vt \in H^1_{tan} (\Gamma)$, we refer by  $\Pi_h \vt$ to a component-wise vector-valued extension of \eqref{globalInt}:
\begin{align}\label{eq:IntComp}
\Pi_h  \vt|_K=\sum^3_{i=1} \phi_i \alpha_i(\tilde \vt),\quad  \alpha_i(\tilde \vt)&=\frac{1}{|E_i|}\int_{E_i} \tilde \vt \, d\sigma_h. 
\end{align}
A decomposition of $\alpha_i(\tilde \vt)$ in normal and conormal components reads as
\begin{align}\label{eq:IntDeComp}
\alpha_i(\tilde \vt)&=\frac{1}{|E_i|}\int_{E_i}  (\tilde{ \vt} \cdot\nt_{E^l}) \nt_{E^l}+ (\tilde{\vt} \cdot\taut_{E^l}) \taut_{E^l}+ (\underbrace{\tilde{ \vt} \cdot \nt}_{=0})\nt \, d\sigma_h.
\end{align}
The vector-valued extension of the scalar interpolation is not necessarily tangential to the discrete surface $\Gamma_h$. In particular $\Pi_h\vt \not\in \Vt_h$.

{
\begin{lemma}[Interpolation Estimates]\label{lemma_INT}
For $\vt \in H^1_{tan}(\Gamma)$ the following interpolation error estimates  hold
\begin{align}\label{eq:H1est}
\| \nabla_{\Gamma_h} (\tilde \vt -{\Pi^{tan}_h \vt} )\|_{L^2(\Gamma_h)}&\leq c h \|\tilde \vt\|_{H^1(\Gamma_h)},\\
\| \Pt(\tilde \vt -\Pi^{tan}_h \vt) \|_{L^2(\Gamma_h)}&\leq c h^2 \|\tilde \vt\|_{H^2(\Gamma_h)},\label{eq:L2est}\\
\|\Pt(\tilde \vt-\Pi^{tan}_h \vt)\|_{L^2(E)}&\leq c h^\frac{3}{2} \|\vt\|_{H^2(\Gamma_h)},
\label{eq:L2estPur}\\
\| \Pt_h(\tilde \vt -\Pi^{tan}_h \vt) \|_{L^2(\Gamma_h)}&\leq c h^2 \|\tilde \vt\|_{H^2(\Gamma_h)},\label{eq:L2estPh}\\
\|\Pt_h(\tilde \vt-\Pi^{tan}_h \vt)\|_{L^2(E)}&\leq c h^\frac{3}{2} \|\vt\|_{H^2(\Gamma_h)},
\label{eq:L2estPurPh}
\end{align}
and in the energy-norm, we have 
\begin{align}\label{eq:int_vec}
\| \tilde \vt-\Pi^{tan}_h  \vt\|_{h}\leq c h \|\tilde \vt\|_{H^2(\Gamma_h)}. 
\end{align}
\end{lemma}

\begin{remark}
To obtain the $L^2$-error estimates with optimal order, we either introduce the projection  $\Pt$ or $\Pt_h$ in front of the difference to account for the fact  that $\tilde \vt=\Pt\tilde \vt$ and $\Pi^{tan}_h \vt=\Pt_h \Pi^{tan}_h \vt$ point in different directions.
\end{remark}
{
\begin{proof}
\textbf{Estimate (\ref{eq:H1est})}: We introduce the interpolation $\Pi_h \vt$.
\begin{equation}
\begin{aligned}
\Big \| &\nabla_{\Gamma_h} (\tilde \vt -\Pi^{tan}_h \vt ) \Big \|_{L^2({K})}\\
&\leq \underbrace{\| \nabla_{\Gamma_h} (\tilde \vt -\Pi_h  \vt )\|_{L^2({K})}}_{O(h),\, (\ref{eq:localerror})} +\| \nabla_{\Gamma_h} (\Pi_h \vt -\Pi^{tan}_h  \vt ) \|_{L^2({K})}.
\end{aligned}
\end{equation}

To bound $\| \nabla_{\Gamma_h} (\Pi_h \vt -\Pi^{tan}_h  \vt ) \|_{L^2({K})}$ we use the expression of the interpolation (\ref{eq:IntComp}), (\ref{eq:IntCompTan})
and use $\nabla_{\Gamma_h} \Pt_h\Pi_h \vt
=\nabla_{\Gamma_h} \Pi_h \vt$, as $\Pt_h$ is constant per cell. Thus,  we get  
\begin{equation}
\begin{aligned}
\Big \|&\nabla_{\Gamma_h} \Big (\Pt_h (\Pi_h \vt -\Pi^{tan}_h \vt)\Big ) \Big  \|_{L^2({K})}\\
&= \Big \|\nabla_{\Gamma_h}\Big(\Pt_h \Big (\sum^3_{i=3} \phi_i (\alpha_i(\tilde \vt)-\alpha_i^{tan}(\tilde \vt) )\Big ) \Big \|_{L^2({K})}\\&=
\Big \|\nabla_{\Gamma_h} \Big (\sum^3_{i=1} \phi_i \Pt_h (\alpha_i (\tilde \vt)-\alpha^{tan}_i(\tilde \vt))\Big)\Big \|_{L^2({K})}\\ &\leq \sum^3_{i=1}{\Big | \Pt_h (\alpha_i(\tilde \vt)-\alpha^{tan}_i(\tilde \vt))\Big|}\cdot{ \Big \|\nabla_{\Gamma_h}  \phi_i \Big \|_{L^2({K})}}\\
 &\leq \sum^3_{i=1} \underbrace{\Big | \Pt_h (\alpha_i(\tilde \vt)-\alpha^{tan}_i(\tilde \vt))\Big|}_{(\star)} \cdot \underbrace{  h^{-1} \Big\|\phi_i \Big \|_{L^2({K})}}_{\leq c\, (\star \star)},
\end{aligned}
\end{equation}
where we used the inverse estimate to reach $(\star \star)$. 
The final estimate results by deriving  a bound for the remaining term $(\star)$. To improve readability, we skip the index $i$ in the following. Using  (\ref{eq:IntComp}) and (\ref{eq:IntCompTan})  we get

\begin{equation}
\begin{aligned}
(\star)&=\Big |\Pt_h \frac{1}{|E|} \int_E (\tilde{ \vt} \cdot\nt_{E^l}) \nt_{E^l}- (\tilde{ \vt} \cdot \nt_{E^l})\nt_E + (\tilde{\vt} \cdot\taut_{E^l}) \taut_{E^l}- (\tilde{\vt} \cdot  \taut_{E^l})\taut_E \, d\sigma_h\Big |\\
&=\Big |\frac{1}{|E|} \int_E (\tilde{ \vt} \cdot\nt_{E^l}) (\Pt_h\tilde \nt_{E^l}-\nt_E) + (\tilde{ \vt} \cdot\taut_{E^l})(\Pt_h\tilde \taut_{E^l}- \taut_E)\, d\sigma_h \Big |\\
&\leq ch^2\frac{1}{h}  \int_{E}|\tilde{ \vt} \cdot \nt_{E^l}| + | \tilde{ \vt} \cdot \taut_{E^l}| \, d\sigma_h 
\leq c h \| \tilde \vt\|_{L^1(E)}\\
&\leq c h^\frac{3}{2} \|\tilde \vt\|_{L^2(E)} 
 \leq c h \|\tilde \vt\|_{H^1(K)},
\end{aligned} 
\end{equation}
where we apply the estimate (\ref{eq:pm_n2}) in the first inequality, the Cauchy-Schwarz inequality in the third inequality, and the trace inequality, as well as Poincare's inequality in the last estimate. 
\textbf{Estimate (\ref{eq:L2est})} is proven similarly. It holds that 
\begin{multline}
\Big \| \Pt(\Pi_h \vt -\Pi^{tan}_h  \vt) \Big  \|_{L^2({K})}= \Big \|\sum^3_{i=1} \phi_i \Pt \Big (\alpha_i (\tilde \vt)-\alpha^{tan}_i(\tilde \vt)\Big)\Big \|_{L^2({K})} \\
\leq \sum^3_{i=1} \Big (\underbrace{\Big | \Pt_h (\alpha_i(\tilde \vt)-\alpha^{tan}_i(\tilde \vt))\Big|}_{\leq ch \|\tilde \vt\|_{H^1(K)}, (\star)}
+ \max_K \underbrace{\|(\Pt-\Pt_h) \|}_{ \leq ch,\, (\ref{eq:P_Ph})} \underbrace{|\alpha_i(\tilde \vt)-\alpha^{tan}_i(\tilde \vt) |}_{ \leq c\|\tilde \vt\|_{H^1(K)}}\Big ) \\ \cdot \underbrace{ \Big \|\phi_i \Big \|_{L^2({K}).}}_{\leq \|1\|_{L^2(K)} \leq 
 c h.}
\end{multline}

\noindent\textbf{Estimate (\ref{eq:L2estPur}):} Invoking the trace inequality, we arrive at
\begin{equation}
\begin{aligned}
\|&\Pt(\tilde \vt -\Pi^{tan}_h \vt)\|_{L^2(E)}\\
&\leq h^{-\frac{1}{2}}\|\Pt(\tilde \vt-\Pi_h^{tan}\vt)\|_{L^2(K)}+h^\frac{1}{2}\|\nabla_{\Gamma_h}\Pt(\tilde \vt-\Pi^{tan}_h \vt)\|_{L^2(K)}.
\end{aligned}
\end{equation}
The first term is bounded by applying (\ref{eq:L2est}). Using the chain rule on the second term and the fact that $\Pt_h\Pt_h=\Pt_h$, we obtain
\begin{equation}
\begin{aligned}
\|\nabla_{\Gamma_h}\Pt (\tilde \vt-\Pi^{tan}_h \vt)\|_{L^2(K)}
&\leq \ |\Pt_h(\Pt -\Pt_h)|_{L^\infty(K)} \| \nabla (\tilde \vt-\Pi^{tan}_h \vt)\Pt_h\|_{L^2(K)}\\
&\quad +\|\Pt_h\nabla(\tilde \vt-\Pi^{tan}_h \vt)\Pt_h\|_{L^2(K)}\\
&\quad + |P|_{W^{1,\infty}} \|\Pt (\tilde \vt-\Pi^{tan}_h \vt)\|_{L^2(K)}\\
& \leq  h \|\tilde \vt\|_{H^2(K)},
\end{aligned} 
\end{equation}
where we use  (\ref{eq:P_Ph}), (\ref{eq:H1est}) and (\ref{eq:L2est}) in the second estimate.
\textbf{Estimate (\ref{eq:L2estPh})} is proven by using the bound $(\star)$.
\textbf{Estimate (\ref{eq:L2estPurPh})} follows by applying (\ref{eq:H1est}).
\textbf{Estimate (\ref{eq:int_vec})} is shown  by combining (\ref{eq:H1est}) and  (\ref{eq:L2est}).
\end{proof}
}
}

\subsection{Energy error estimate}\label{sec;energy}
The aim of this section is to derive the energy error in the discrete energy norm. The main tool used in the derivation is Strang's second Lemma. 
\begin{theorem}[Strang's second Lemma]\label{eq:strang}
Let  $\ut\in H^1_{tan}(\Gamma)$ be the exact solution of (\ref{eq:ana}), $\tilde \ut$ be its extension to $U$ and $\ut_h \in \Vt_h$ be the discrete solution of (\ref{eq:var}). Then the following estimate holds 
\begin{equation}
  \|\tilde \ut-\ut_h\|_h \leq \underbrace{c \inf_{\vt_h \in \Vt_h} \|\tilde \ut-\vt_h\|_h}_{\text{approximation error}} + \underbrace{\sup_{\vt_h \in \Vt_h} \frac{a_h(\tilde \ut, \vt_h) - (\tilde \ft,\vt_h)}{\| \vt_h\|_h}}_{\text{nonconformity error}}.
\end{equation}
\end{theorem}
 The following lemma proves the nonconformity error. 
 
  \begin{lemma}[Nonconformity Error Estimate]\label{lemma:noncon}
 Let $\ut$ be the solution of (\ref{eq:ana}) and $ \tilde \ut$ be its extension to $U$ defined by (\ref{eq:extension}). Then the following error estimate holds
 \begin{equation}\label{eq:noncon}
    |a_h(\tilde \ut,\vt_h)-(\tilde \ft,\vt_h)_{\Gamma_h}| \leq c h \|\ft\|_{L^2(\Gamma)}\|\vt_h\|_h,
 \end{equation}
 for any $\vt_h \in \Vt_h$.
 \end{lemma}
 { \begin{proof} By inserting $(\ft,\vt_h^ l)_\Gamma$  we get 
  \begin{align}\label{eq:first}
    |a_h(\tilde \ut,\vt_h)-(\tilde \ft,\vt_h)_{\Gamma_h}|=& |a_h(\tilde \ut,\vt_h)-(\ft,\vt_h^l)_\Gamma + (\ft,\vt_h^l)_\Gamma-(\tilde \ft,\vt_h)_{\Gamma_h}|.
    \end{align}
     {
We use Green's first identity 
\begin{equation}\label{eq:green}
  { \Big ( \nabla_{\Gamma}\ut,\nabla_{\Gamma}\vt^ l_h \Big)_{K^l}= -\Big (\Pt \div_{\Gamma}(\nabla_{\Gamma} \ut), \vt^l_h \Big)_{K^l}
 +\int_{\partial K^l} (\nabla_{\Gamma} \ut) \nt_{E^l} \cdot \vt^l_h} d \sigma,
\end{equation}}
 and obtain 
    \begin{align}\label{eq:first}
    |a_h(\tilde \ut,\vt_h)-(\tilde \ft,\vt_h)_{\Gamma_h}|
    \leq & \underbrace{\Big |a_h(\tilde \ut,\vt_h)- a(\ut,\vt_h^l)+ (\ft,\vt_h^l)_{\Gamma}-(\tilde \ft,\vt_h)_{\Gamma_h} \Big|}_{I}\\
    &+\underbrace{\sum_{E^l}\Big |\int_{E^l}\vt_h^{+,l} \cdot \nabla_\Gamma \ut \nt^+_{E^l} + \vt_h^{-,l} \cdot \nabla_\Gamma \ut \nt^-_{E^l} \, d\sigma\Big|}_{II},
    \end{align} 
    where $\vt_h^{+,l}:=\vt^l_h|_{K^+}$.
    The terms collected in $I$ are the geometric error, whereas the jump terms that arise due to the nonconformity of the element are collected in $II$.\\
   \textbf{Geometric error $I$:}
The estimate for $I$ is derived by applying the estimates provided in Lemma \ref{Lemma:Geo}, the norm equivalence given in Lemma \ref{lemma:norm}, and the regularity estimate (\ref{eq:regu}):
   \begin{equation}\label{eq:geo_Q}
  I=  |a_h(\tilde \ut,\vt_h)- a(\ut,\vt_h^l)+ (\ft,\vt_h^l)_\Gamma-(\tilde \ft,\vt_h)_{\Gamma_h}|
    \leq h \|\vt_h\|_h\|\ft\|_{L^2(\Gamma)}.
    \end{equation}
  \textbf{Jump terms $II$:} The terms can  be  bounded by 
 \begin{equation}
\begin{aligned}
\sum_{E^l}&\Big |\int_{E^l}[\vt_h^{l}] \cdot \nabla_\Gamma \ut \nt^+_{E^l} d \sigma \Big |\\
&\leq \underbrace{\sum_{E^l}\Big |\int_{E^l}[v_n^{l}\nt^l_E] \cdot  \nabla_\Gamma \ut \nt^+_{E^l} d\sigma\Big|}_{II_1}+ \underbrace{\sum_{E^l}\Big|\int_{E^l}[v_t^l \taut^l_E] \cdot  \nabla_\Gamma \ut \nt^+_{E^l} d \sigma \Big |}_{II_2},
\end{aligned}
\end{equation}
where 
$\int_{E^l}[v_n\nt_E] d\sigma =\int_{E^l}v^{+,l}_n
\nt^{+,l}_E-v^{-,l}_n\nt^{-,l}_E d\sigma$ 
and  $\int_{E^l}[v_n\taut_E] d\sigma=\int_{E^l}v^{+,l}_n\taut^{+,l}_E-v^{-,l}_n\taut^{-,l}_E d\sigma$.
Note that $\nt_{E^l}$ is the  outer conormal vector to $K^l$ on $\Gamma$,
while $\nt^l_E$ is the lift of $\nt_E$ from $\Gamma_h$ 
to $\Gamma$.  It holds that $\nt^+_{E^l}=-\nt^-_{E^l}$, but $\nt^+_E\neq - 
\nt^-_E$.  \\
\textbf{Bound $II_1$}.  We use the fact that $\Pt^T=\Pt$ and $\Pt \nabla_\Gamma \ut=\nabla_\Gamma \ut$.
 We add and subtract suitable terms and get  
 \begin{equation}\label{eq:II1}
 \begin{aligned}
II_1&=\sum_{E^l}\Big|\int_{E^l} \Big ( \Pt( v^{+,l}_n\nt^{+,l}_E)- \Pt(v^{-,l}_n\nt^{-,l}_E) \Big ) \cdot \nabla_\Gamma \ut \nt^+_{E^l} \, d \sigma \Big|\\
 & \leq \sum_{E^l}\Big| \int_{E^l} \Big ( v^{+,l}_n (\Pt \nt^{+,l}_E- \nt^+_{E^l}) - v^{-,l}_n(\Pt \nt^{-,l}_E -\nt^ -_{E^l })  \Big ) \cdot \nabla_\Gamma\ut \nt^+_{E^l} \, d \sigma \Big|\\
& + \sum_{E^l}\Big |\int_{E^l}(v^{+,l}_n\nt^+_{E^l}-v^{-,l}_n\nt^-_{E^l})\cdot \nabla_\Gamma \ut \nt^+_{E^l} d \sigma  \Big|\\
&\leq \sum_{E^l} \underbrace{\|\Pt \nt^{+,l}_E- \nt^+_{E^l}\|_{L^\infty(E^l)}}_{O(h^2),\, (\ref{eq:pm_n}) }\|v^{+,l}_n\|_{L^2(E^l)} \| \nabla_\Gamma\ut\|_{L^2(E^l)}\ \\
&+\sum_{E^l} \underbrace{\|\Pt \nt^{-,l}_E -\nt^ -_{E^l }\|_{L^\infty(E^l)}}_{O(h^2), \, (\ref{eq:pm_n})}{\|v^{-,l}_n\|_{L^2(E^l)} \| \nabla_\Gamma\ut\|_{L^2(E^l)}}\\
&+ \sum_{E^l}\Big |{\int_{E^l}{[v_n^l]} \nt^+_{E^l} \cdot \nabla_\Gamma\ut \nt^+_{E^l} \, d \sigma }\Big|\\
&\leq c  h \| \vt_h\|_{h}\|\ft\|_{L^2(\Gamma)},
 \end{aligned}
 \end{equation}
where we apply {the trace inequality and the regularity estimate (\ref{eq:regu}). To bound
\[
\sum_{E^l}\Big |{\int_{E^l}{[v_n^l]} \nt^+_{E^l} \cdot \nabla_\Gamma\ut \nt^+_{E^l} \, d \sigma }\Big| \leq c  h \| \vt_h\|_{h}\|\ft\|_{L^2(\Gamma)}
\]
we refer to Lemma \ref{lemma:jumpSmall} which we moved to the appendix due to its technical nature.} 
\textbf{A bound for $II_2$} is derived in a similar way. 
\end{proof}

 We now turn to presenting the main error estimate of Section \ref{sec;energy}, the energy error estimate. 
  
   \begin{theorem}[Energy-error estimate]\label{th:energy} Let $\tilde \ut$ be the extension of  the exact solution of (\ref{eq:ana}) to U and $\ut_h$ the discrete solution of (\ref{eq:var}). Then the following error estimate holds
  \begin{equation}\label{eq_energy}
  \|\tilde \ut-\ut_h\|_h \leq c h \|\ft\|_{L^2(\Gamma)}.
  \end{equation}
  \end{theorem}
  \begin{proof}
The nonconformity estimate provided in Lemma \ref{lemma:noncon}  gives 
  \begin{equation}\label{eq:th}
  \sup_{\vt_h \in \Vt_h} \frac{|a_h(\tilde \ut, \vt_h)- (\ft,\vt^l_h)|}{\|\vt_h\|_h} \leq c h \|\ft \|_{L^2(\Gamma)}. 
  \end{equation}
  
Next, we consider the approximation error:
{\Gr
\begin{equation}
\begin{aligned}\label{eq:approxi}
    \inf_{\vt_h \in \Vt_h}\| \tilde \ut-\vt_h\|_h& \leq \| \tilde \ut-\Pi^{tan}_h  \ut \|_h\leq c h\|\tilde \ut\|_{H^2(\Gamma_h)}\leq c h\|\ft\|_{L^ 2(\Gamma)},
\end{aligned}
\end{equation}
}
where we used the interpolation estimate (\ref{eq:int_vec}) in the second bound and the regularity estimate   (\ref{eq:regu}) in the last inequality.
The proof is completed by the combination of Strang's second Lemma (\ref{eq:strang}) with the error bounds  (\ref{eq:th}) and (\ref{eq:approxi}).
  \end{proof}
   \subsection{$L^2$-error estimate}\label{Sec:L2}
  The derivation of the $L^2$-error estimate is based on the Aubin-Nitsche trick. We introduce the dual problem 
  \begin{equation}\label{eq:ana_dual}
  \begin{aligned}
   \text{Find } \zt \in H_{tan}^1(\Gamma) \text{ s.t. } \quad a(\vt,\zt)& =(\vt,\gt),\quad \forall \vt \in H^1_{tan}(\Gamma),\\ \gt&={\Pt_h}(\ut-\ut_h^l) \in L^2(\Gamma).    
  \end{aligned}
  \end{equation}
%
The following regularity estimate holds (cf. \cite[eq. (5.158)]{HansboLarson2020}).
  \begin{equation}\label{eq:reg_dual}
  \|\zt\|_{H^2({\Gamma})} \leq  c \|\gt\|_{L^2(\Gamma)}.
  \end{equation}
The corresponding discrete dual problem  to (\ref{eq:ana_dual}) is defined as 
  \begin{equation}\label{eq:dual_disc}
\text{find } \zt_h \in \Vt_h \text{ s.t.}\quad  a_h(\vt_h,\zt_h)=(\vt_h, \tilde \gt)_{\Gamma_h}, \quad \forall \vt_h \in \Vt_h.
\end{equation}  
Based on the energy-error estimate provided in Theorem \ref{th:energy} the following energy-estimate results for the dual solution 
 \begin{equation}\label{eq:dual_energy}
 \| \tilde \zt-\zt_h\|_h \leq c h \|\gt\|_{L^2(\Gamma)}.
 \end{equation}
 {
We present the following primal and dual interpolation estimate involving the interpolation $\Pi_h^{tan}$.
 \begin{lemma}[Primal and dual interpolation estimate]\label{lemma:Klaus2}
 Let $\ut$ be the solution of the primal problem (\ref{eq:ana}) and $\zt \in H^2(\Gamma) $ be a vector field. It holds that
\begin{align}
a_h(\tilde \ut,\tilde \zt-{\Pi^{tan}_h} \zt) &\leq c h^2 \|\ft\|_{L^2(\Gamma)} \|\zt\|_{H^2(\Gamma)}.
\end{align}
The following estimate holds for a vector field $\ut \in H^2(\Gamma)$ and the dual solution $\zt$ of (\ref{eq:reg_dual})
\begin{align}
a_h(\tilde \ut-{\Pi^{tan}_h} \ut,\tilde \zt) &\leq c h^2 \|\ut\|_{H^2(\Gamma)} \|\gt\|_{L^2(\Gamma)}.
\end{align}
 \end{lemma}
  \begin{proof}
  We prove this estimate in Section  \ref{sec:appbig} of the Appendix.
 \end{proof}

Based on the previous lemma, we next show estimates for the consistency error.
{
\begin{lemma}[Primal and dual consistency error estimate]\label{lemmma:dual_consist}
    Let $\ut$ be the solution to the primal problem (\ref{eq:ana}) and  $\zt$
     be the solution to the dual problem $(\ref{eq:ana_dual})$. $\ut_h$ is the discrete solution of (\ref{eq:var}) and $\zt_h$ the discrete solution of (\ref{eq:dual_disc}). Then the following estimates hold
    \begin{align*}
       \Big | a_h(\tilde \ut, \tilde \zt-\zt_h) -\Big ((\ft,\zt)_\Gamma- (\tilde \ft ,\zt_h)_{\Gamma_h}\Big ) \Big| &\leq c h^2 \|\ft\|_{L^2(\Gamma)}\|\gt\| _{L^2(\Gamma)},\\
       \Big | a_h(\tilde \ut-\ut_h,\tilde \zt) -\Big( (\ut,\gt)_\Gamma-(\ut_h,\tilde \gt)_{\Gamma_h} \Big) \Big | &\leq c h^2 \|\ft\|_{L^2(\Gamma)}\|\gt\|_{L^2(\Gamma)}.\\
    \end{align*}
\end{lemma}
 \begin{proof}
     Both inequalities can be proven analogously. We outline the proof of the first one and introduce the interpolation $\Pi_h^{tan}\zt$ as well as $\tilde \zt$. 
     \begin{align*}
      a_h(\tilde \ut, &\tilde \zt-\zt_h) -\Big ((\ft,\zt)_\Gamma- (\tilde \ft ,\zt_h)_{\Gamma_h}\Big )\\
   &=a_h(\tilde \ut, \tilde \zt- \Pi_h^{tan} \zt)+a_h(\tilde \ut,  \Pi_h^{tan}\zt -\zt_h) - (\ft,\zt)_\Gamma+ (\tilde \ft ,\zt_h -\Pi_h^{tan}\zt)_{\Gamma_h}\\
   &\quad + (\tilde \ft ,\Pi_h^{tan}\zt -\tilde \zt)_{\Gamma_h} +(\tilde \ft ,\tilde \zt)_{\Gamma_h}\\           
         &=\underbrace{a_h(\tilde \ut, \tilde \zt-{\Pi^{tan}_h} \zt)}_{I_1}
        +\underbrace{(\tilde \ft,\tilde \zt)_{\Gamma_h}-(\ft,\zt)_\Gamma}_{I_2}+\underbrace{(\tilde \ft,-\tilde \zt+ {\Pi^{tan}_h} \zt)_{\Gamma_h}}_{I_3}\\
         &\quad  +\underbrace{a_h( \tilde \ut,{\Pi^{tan}_h} \zt-\zt_h)+( \tilde \ft,\zt_h-{\Pi^{tan}_h}\zt)_{\Gamma_h}}_{I_4}.
     \end{align*}
Using the primal interpolation estimate from Lemma \ref{lemma:Klaus2}, $I_1$ is bound by
\begin{align*}
    |I_1|=a_h(\tilde \ut, \tilde \zt- {\Pi^{tan}_h}\zt)\leq c h^2 \|\ft\|_{L^2(\Gamma)} \|\zt\|_{H^2(\Gamma)}.
\end{align*}
The  geometric estimates provided in Lemma \ref{Lemma:Geo} give
  \begin{align*}
      |I_2|=|(\tilde \ft, \tilde \zt)_{\Gamma_h}-(\ft,\zt)_\Gamma|\leq c h^2 \|\tilde \ft\|_{L^2(\Gamma_h)}\|\tilde \zt\|_{L^2(\Gamma_h)}\leq c h^2 \|\tilde \ft\|_{L^2(\Gamma_h)}\|\tilde \zt\|_{H^2(\Gamma_h)}.
  \end{align*}
Combining the Cauchy-Schwarz inequality, the interpolation estimate (\ref{eq:L2est}) and  the fact that $\tilde \ft=\Pt \tilde \ft$, $\tilde \zt=\Pt \tilde \zt$, $\Pt^T=\Pt$  gives 
\begin{align*}
    |I_3|=\Big|(\Pt\tilde \ft,{\Pi^{tan}_h}\zt- \tilde \zt)_{\Gamma_h}\Big |&\leq \Big (\int_{\Gamma_h} \tilde \ft^2 d \sigma_h\Big )^\frac{1}{2}\Big(\int_{\Gamma_h} |\Pt({\Pi^{tan}_h}\zt - \tilde \zt)|^2 d \sigma_h \Big)^\frac{1}{2}\\
    &\leq c h^2 \| \tilde \ft \|_{L^2(\Gamma)} \|\tilde \zt\|_{H^2(\Gamma_h)}.
\end{align*}
We bound $I_4$ by applying the nonconformity estimate given in Lemma \ref{lemma:noncon}  
  \begin{equation*}
      |I_4|=|a_h( \tilde \ut,{\Pi^{tan}_h} \zt-\zt_h)+( \tilde \ft,\zt_h-{\Pi^{tan}_h}\zt)_{\Gamma_h}|
      \leq c h \|\ft \|_{L^2(\Gamma)}\|{\Pi^{tan}_h} \zt -\zt_h \|_h . 
  \end{equation*}
  Adding and subtracting suitable terms, as well as applying the  interpolation  estimate in the energy norm (\ref{eq:int_vec}) and  the energy estimate (\ref{eq:dual_energy}) yields
   \begin{align*}
      |I_4|&\leq c h \|\ft \|_{L^2(\Gamma)} \Big( \|{\Pi^{tan}_h} \zt- \tilde \zt \|_h+\| \tilde \zt -\zt_h  \|_h\Big)\\
      &\leq c h \|\ft \|_{L^2(\Gamma_h)} {\Big( ch \|\tilde \zt\|_{H^2(\Gamma_h)}+h\|\gt\|_{L^2(\Gamma)}\Big)}.
  \end{align*}
 The final estimate results by applying the norm equivalence provided in Lemma \ref{lemma:norm} and the regularity estimate (\ref{eq:reg_dual}) to the bounds derived for $I_1$ to $I_4$.
 \end{proof}
}
\begin{theorem}[$L^2$ Error Estimate]\label{th:L2}
    Let $\ut$ be the primal solution to (\ref{eq:ana}), $ \tilde \ut$ be the extension to a neighborhood and $\ut_h$ be the discrete solution of (\ref{eq:var}). Then the following error estimate holds 
    \begin{align*}
        \| \Pt_h (\tilde \ut -\ut_h)\|_{L^2(\Gamma_h)} \leq c h^2 \|\ft\|_{L^2(\Gamma)}.
    \end{align*}
\end{theorem}

\begin{proof}
Note that $\Pt_h \tilde \gt=\tilde \gt$. First, we introduce the geometrical error, collected in terms $I_1,I_2$, by adding and subtracting suitable terms. 
    \begin{align*}
        \| \Pt_h(\tilde \ut-\ut_h)\|^2_{L^2(\Gamma_h)}=&(\tilde \ut-\ut_h, \Pt_h \tilde \gt)_{\Gamma_h}\\
        =&\underbrace{(\tilde \ut, \tilde \gt)_{\Gamma_h} - (\ut,\gt)_\Gamma}_{I_1} + (\ut,\gt)_\Gamma-(\ut_h, \tilde \gt )_{\Gamma_h}\underbrace{-a(\ut,\zt) +a_h( \tilde \ut, \tilde \zt)}_{I_2}\\
        &-a_h(\tilde \ut,\tilde \zt) +a(\ut,\zt).
        \end{align*}
        In the second step, we add and subtract further terms to apply  the primal and dual energy-error estimate to the resulting term $I_3$.
         \begin{align*}
        \| \Pt_h(\tilde \ut-\ut_h)\|^2_{L^2(\Gamma_h)}
        =&I_1+I_2+(\ut,\gt)_\Gamma-(\ut_h,\tilde \gt )_{\Gamma_h}-a_h( \tilde \ut, \tilde \zt)+a(\ut,\zt)\\
        &+ \underbrace{a_h(\tilde \ut-\ut_h, \tilde \zt - \zt_h)}_{I_3}-a_h(\tilde \ut, \tilde \zt)+ a_h(\ut_h,\tilde \zt)\\&+a_h(\tilde \ut,\zt_h)-a_h(\ut_h,\zt_h).
\end{align*}         
       We make use of  the relation of (\ref{eq:ana}) and (\ref{eq:ana_disc}) and reorder the terms to insert the primal and dual consistency error, collected in term $I_4$ and term $I_5$, respectively.
\begin{align*}
  \| \Pt_h(\tilde \ut-\ut_h)\|^2_{L^2(\Gamma_h)}       &=I_1+I_2+I_3+(\ut,\gt)_\Gamma-(\ut_h, \tilde \gt )_{\Gamma_h}-a_h(\tilde \ut, \tilde \zt)+(\ft,\zt)_\Gamma\\
        &\quad -a_h(\tilde \ut, \tilde \zt)+ a_h(\ut_h,\tilde \zt)+a_h(\tilde \ut,\zt_h)-(\tilde \ft,\zt_h)_{\Gamma_h}\\
        &= I_1+I_2+I_3\underbrace{-a_h(\tilde \ut-\ut_h,\tilde \zt)+(\ut,\gt)_\Gamma-(\ut_h,\tilde \gt )_{\Gamma_h}}_{I_5}\\
        &\quad \underbrace{-a_h(\tilde \ut,\tilde \zt -\zt_h)+(\ft,\zt)_\Gamma-(\tilde \ft,\zt_h)_{\Gamma_h}}_{I_4}.
    \end{align*}    
Based on the geometric estimates provided in Lemma \ref{Lemma:Geo}  and by  applying the norm equivalences provided in Lemma \ref{lemma:norm} as well as the regularity estimates (\ref{eq:regu}) and (\ref{eq:reg_dual})  one gets 
    \begin{align}
        |I_1|&\leq c h^2 \|\tilde \ut\|_{L^2(\Gamma_h)}\|\tilde \gt\|_{L^2(\Gamma_h)}\leq c h^2 \|\ft\|_{L^2(\Gamma)}\|\gt\|_{L^2(\Gamma)},\\
         |I_2|&\leq c h^2 \|\ut\|_{H^ 2(\Gamma)}\|\zt\|_{H^2(\Gamma)} \leq c h^ 2\|\ft\|_{L^2(\Gamma)}\|\gt\|_{L^ 2(\Gamma)}.
    \end{align}
     We bound $I_3$ by applying the primal (\ref{eq_energy}) and dual energy estimate (\ref{eq:dual_energy}) and obtain
\begin{align*}
    |I_3|
    \leq c& \|\tilde \ut-\ut_h\|_h \|\tilde \zt-\zt_h\|_h
    \leq c h^2 \|\ft\|_{L^2(\Gamma)}\|\gt\|_{L^2(\Gamma)}.
\end{align*}
We apply the bounds for the primal and dual consistency error provided in Lemma \ref{lemmma:dual_consist} 
 to $I_4$ and $I_5$ and get 
\begin{align*}
    |I_4|+|I_5|\leq  c h^ 2 \|\ft\|_{L^ 2(\Gamma)}\|\gt\|_{L^2(\Gamma)}.
\end{align*}
The combination of the bounds for the terms $I_1$ to $I_5$ gives the final estimate.
\end{proof}
\section{Numerical results}\label{sec:num}
As illustrative cases, we consider the surface of the unit sphere and a torus with an inner radius one and an outer radius one half, both representing two-dimensional manifolds embedded in $\mathbb{R}^3$ and characterized by their implicit descriptions $\Gamma_S$ and $\Gamma_T$, respectively.
\begin{align*}
\Gamma_{S}&=\Big \{ \xt=(x_1,x_2,x_3)^T \in \mathbb{R}^3|\,\phit_S(\xt):=(x_1^2+x_2^2+x_3^2-1)=0 \Big \},\\
\Gamma_{T}&=\Big \{ \xt=(x_1,x_2,x_3)^T \in \mathbb{R}^3|\,\phit_T(\xt):=-\Big (\frac{1}{4} - x_3^2 - (\sqrt{x_1^2+x_2^2}-1)^2\Big )=0 \Big \}.
\end{align*}
This level-set formulation enables an analytical computation of the normal vector $\nt$.

\begin{figure}
\begin{center}
\begin{tabular}{rcrc}
a)&&b)\\
&\includegraphics[width=0.4\textwidth]{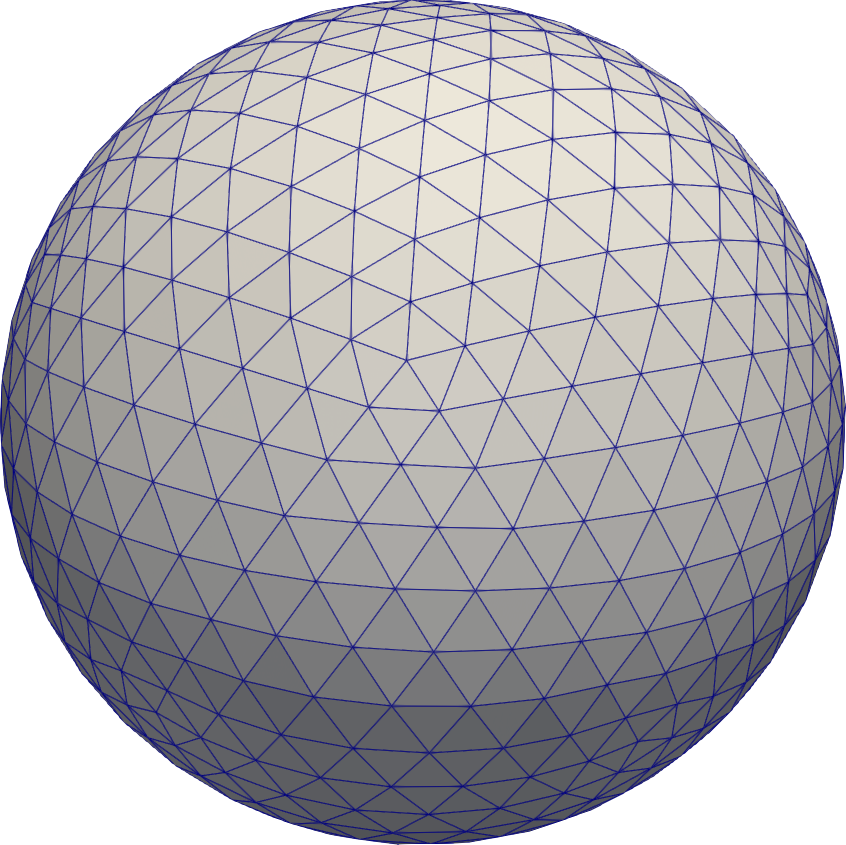}&&\includegraphics[width=0.4\textwidth]{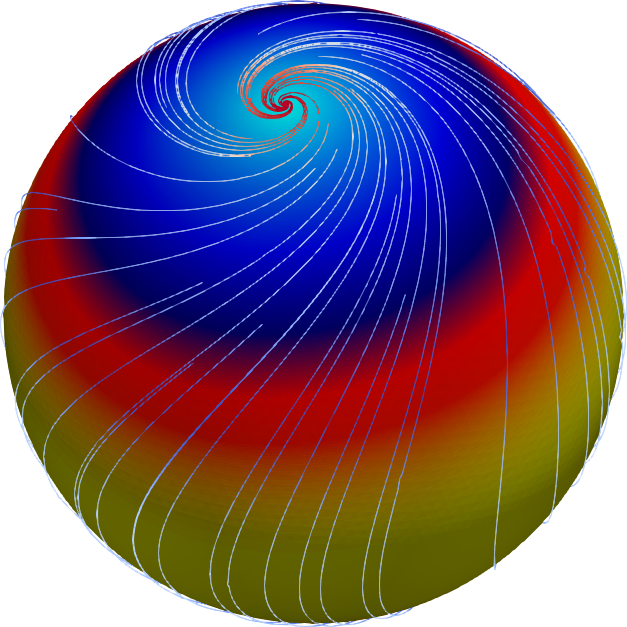}
\end{tabular}
\end{center}
\caption{Sphere: a): Grid after three refinements. b): First component of the solution $\ut_{\Gamma_S}$ with vorticity streamlines.\label{fig:sphere}}
\end{figure}
We implement the surface Crouzeix-Raviart element in Python, following the same implementation as in the climate model ICON \cite{MehlmannGutjahr2022}. In contrast to the latter, which is restricted to approximating the Earth by a sphere, the Python realization offers a larger flexibility and can be applied to a variety of manifolds such as a torus. 
The discrete surface approximation $\Gamma_h$ is generated using the Python package \emph{PyVista} \cite{pyvista}. To evaluate the errors, we prescribe an exact solution of (\ref{eq:conEq}) and determine the corresponding right-hand side function $f$ analytically using \emph{SymPy} \cite{Sympy}. The solutions $\ut_{\Gamma_S}$ and $\ut_{\Gamma_T}$ are prescribed on the surfaces of the sphere and the torus, respectively, as
\begin{align}
\ut_{\Gamma_S}(\xt)&=\Pt(\sin(x_2x_3),-\sin(x_1x_3),\cos(x_3^2))^T,\quad \text{for } \xt \in \Gamma_{S},\label{eq:sphere}\\
\ut_{\Gamma_T}(\xt)&=\Pt(x_2+x_3x_1,-x_1x_3,x_3^2)^T,\quad \text{for } \xt \in \Gamma_{T}\label{eq:torus}.
\end{align}
The right-hand side of (\ref{eq:conEq}) is given as  $f_{\Gamma_S}=-\Delta_\Gamma\ut_{\Gamma_S} +\frac{1}{10}\ut_{\Gamma_S}$ and  $f_{\Gamma_T}=-\Delta_\Gamma\ut_{\Gamma_T} +\frac{1}{10}\ut_{\Gamma_T}$, respectively. In Figure \ref{fig:sphere} and Figure \ref{fig:torus} the
surface geometry and solution on the sphere and the torus is visualized, respectively. 
\begin{figure}
\begin{tabular}{l l}
a)&b)\\
\includegraphics[width=0.43\textwidth]{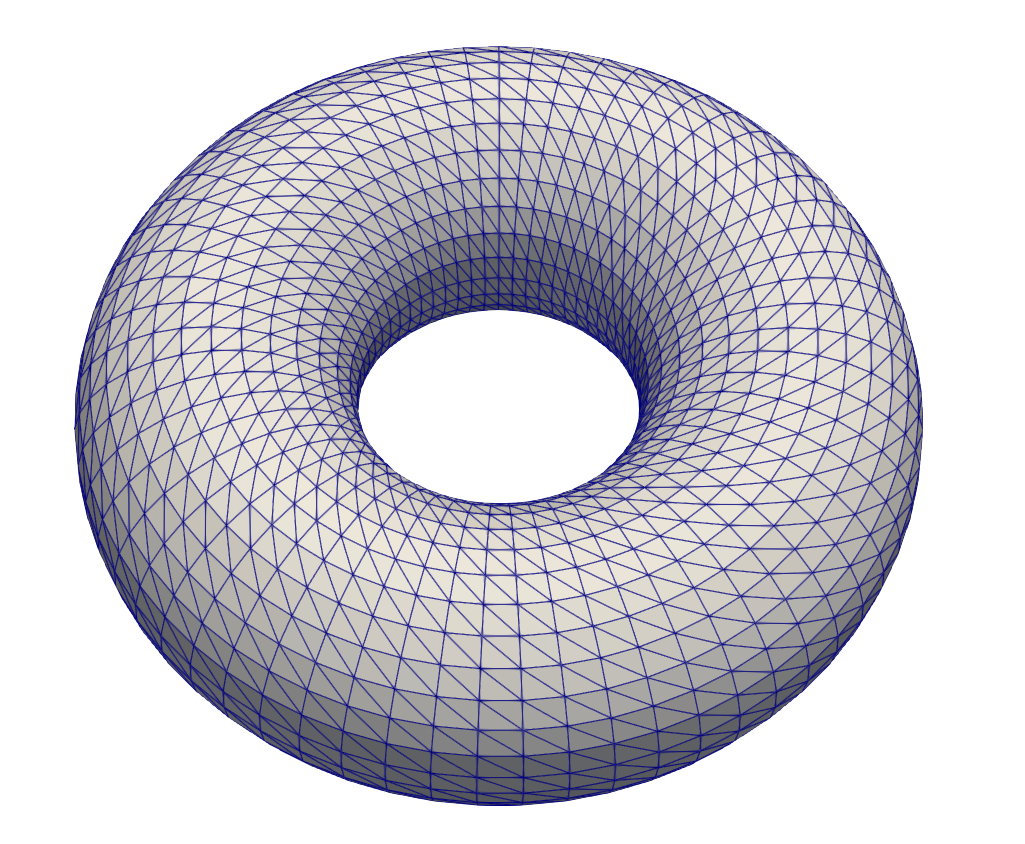}&\includegraphics[width=0.45\textwidth]{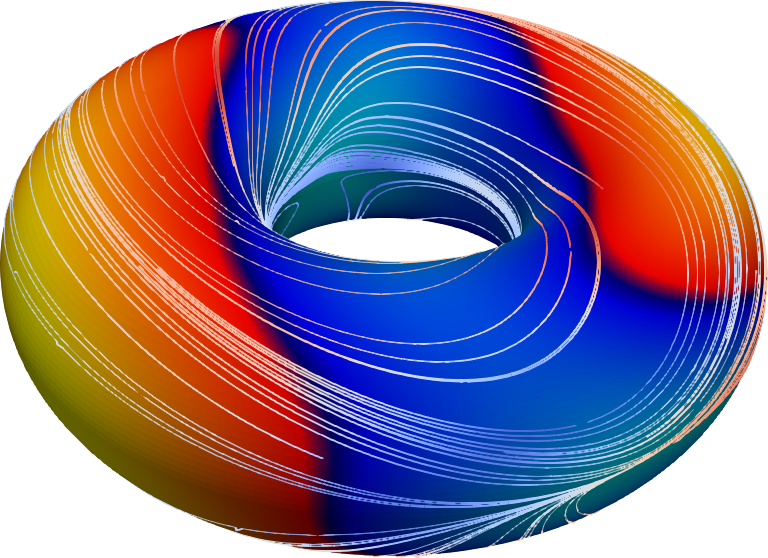}
\end{tabular}
\caption{Torus: a): Grid after three refinements. b): First component of the solution $\ut_{\Gamma_T}$ with vorticity streamlines.\label{fig:torus}}
\end{figure}

\subsection{Numerical results on the sphere}
 The unit sphere is approximated by triangular elements in a sequence of meshes ranging from $h=0.325$ to $h=0.021$. 
The error between the extended exact and the approximate  solution, $\tilde \ut -\ut_h$, is evaluated in the $H^1$-norm. According to Theorem \ref{th:L2} we consider $\Pt_h(\tilde \ut-\ut_h)$ in the $L^2$-norm. We present the numerical results in Table \ref{tab:sphere}. In line with Theorem \ref{th:L2} and Theorem \ref{th:energy} the L2-error converges with quadratic order, $O(h^2),$ and the $H^1$-norm converges with first order, $O(h)$. The errors under mesh refinement are visualized in the left panel of Figure \ref{convOrder}.
\begin{table}[htbp]
\caption{Sphere: Experimental order of convergence on the unit sphere  for the analytic solution (\ref{eq:sphere}).\label{tab:sphere}}
\begin{center}
  \begin{tabular}{r |c c |c c}
  \toprule
   Edge length& $L^2$-error& Order &$H^1$-error& Order\\ 
\midrule
0.325 & $5.66\cdot 10^{-2}$ & --- & $3.39\cdot 10^{-1}$ & ---\\
0.165 & $1.40\cdot 10^{-2}$ & 2.02 & $1.62\cdot 10^{-2}$ & 1.06\\
0.083 & $3.49\cdot 10^{-3}$ & 2.00 & $7.99\cdot 10^{-2}$& 1.02\\
0.041 & $8.71\cdot 10^{-4}$ & 2.00 & $3.98\cdot 10^{-2}$ & 1.01\\
0.021 & $2.20\cdot 10^{-4}$ & 1.99 & $1.99\cdot 10^{-2}$ & 1.00\\
\bottomrule
  \end{tabular}
\end{center}
\end{table}

\subsection{Numerical results on the torus}
 As on the sphere, we consider triangular elements with side lengths between $0.364$ and $0.022$ to discretize the torus.
The error between the exact and numerical solutions is assessed in the $H^1$-norm and in the $L^2$-norm,  where we also evaluate $\Pt_h(\tilde \ut - \ut_h)$. Table \ref{tab:torus} summarizes the numerical results. According  to Theorem \ref{th:energy} and Theorem \ref{th:L2} the $H^1$-error converges linearly and the $L^2$-error converges quadratically. The right panel of Figure \ref{convOrder} shows the error under mesh refinement.  

\begin{table}[htbp]
\caption{Torus: Experimental order of convergence on the torus  for the analytic solution (\ref{eq:torus}).\label{tab:torus}}
\begin{center}
  \begin{tabular}{r |c c |c c}
  \toprule
   Edge length& $L^2$-error& Order &$H^1$-error& Order\\ 
\midrule
0.364 &  $ 1.11\cdot 10^{-1}$ & --- & $1.00 \cdot 10^{-0}$ & ---\\
0.180 & $3.07\cdot 10^{-2}$ & 1.85 & $4.76\cdot 10^{-1}$ & 1.08\\
0.089 & $7.77\cdot 10^{-3}$& 1.98 & $2.33\cdot 10^{-1}$ & 1.03\\
0.044 & $1.94\cdot 10^{-3}$ & 2.00 & $1.15\cdot 10^{-1}$ & 1.01\\
0.022 & $4.84\cdot 10^{-4}$ & 2.00 & $5.75\cdot 10^{-1}$ & 1.01\\
\bottomrule
  \end{tabular}
\end{center}
\end{table}

\begin{figure}
\begin{tabular}{l l}
a)&b)\\
\includegraphics[width=0.46\textwidth]{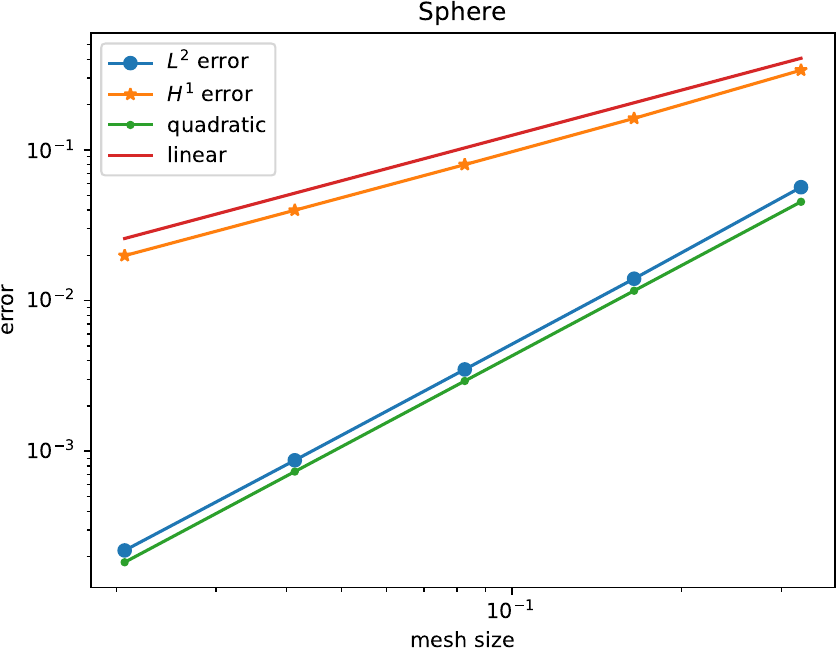}&\includegraphics[width=0.46\textwidth]{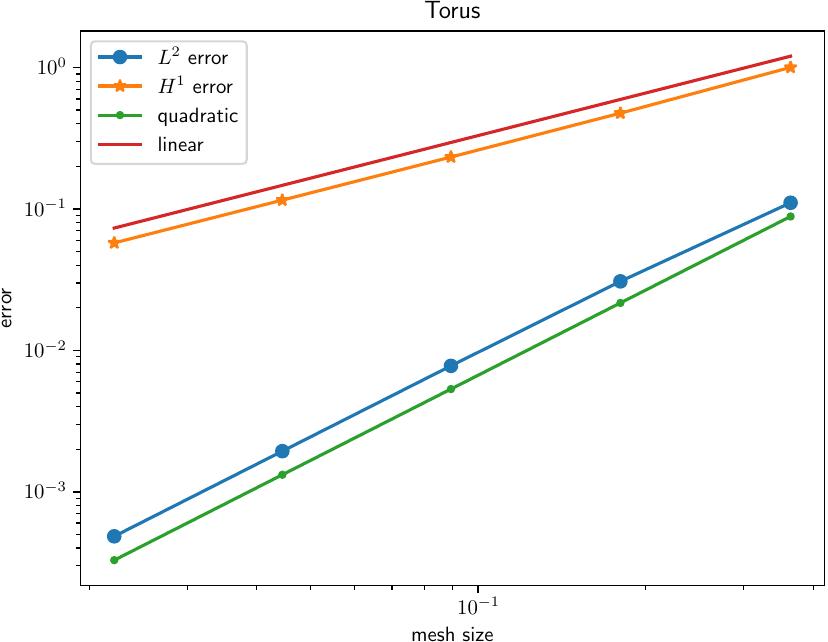}
\end{tabular}
\caption{Errors in the norms $\| \tilde \ut -\ut_h\|_{H^1(\Gamma_h)} $ and $ \|\Pt_h(\tilde \ut -\ut_h)\|_{L^2(\Gamma_h)}$ on the sphere and on the tours. The absolute values are given in Table \ref{tab:sphere} and Table \ref{tab:torus}, respectively. \label{convOrder}}
\end{figure}



\appendix
\section{Appendix}\label{sec:app}
\subsection{ Estimates of geometric errors}
\begin{lemma}[Operator differences]\label{lemma:OpDiff}
For $ \vt \in H^1(\Gamma_h)$  we get 
\begin{equation}
\| (\nabla_\Gamma \vt^l)^e-\nabla_{\Gamma_h}\vt\|_{L^2(\Gamma_h)}\leq h \|\vt\|_{H^1(\Gamma_h)},
\end{equation}
where $e$ indicates the extension.
\end{lemma}
\begin{proof}
The proof follows from \cite{HansboLarson2020} Lemma 5.4 with $k_g=1$ and $\vt_{n_h}=0$.
\end{proof}

\begin{lemma}[Geometric estimates]\label{Lemma:Geo}
 For $\ut, \vt \in H^1(\Gamma_h)$ and $ \ft, \wt \in L^2(\Gamma_h)$ the following estimates hold
 \begin{align}
   | (\ft^l,\wt^l)_\Gamma-(\ft,\wt)_{\Gamma_h}|
    \leq& h^2 \|\wt\|_{L^2(\Gamma_h)}\|\ft\|_{L^2(\Gamma_h)},\label{eq:geomF}\\
    |(\nabla_{\Gamma_h}\ut,\nabla_{\Gamma_h} \vt)_{\Gamma_h}- (\nabla_\Gamma\ut,\nabla_\Gamma \vt_h^l)_\Gamma|
    \leq& {h} \|\vt\|_{H^1(\Gamma_h)}\|\ut\|_{H^1(\Gamma_h)}.\label{eq:geomA}
 \end{align}

 and for $\ut,\vt \in H^2_{tan}(\Gamma)$ it holds that 

 \begin{align}
    |(\nabla_{\Gamma_h}\ut,\nabla_{\Gamma_h} \vt)_{\Gamma_h}- (\nabla_\Gamma\ut,\nabla_\Gamma \vt_h^l)_\Gamma|
    \leq& {h^2} \|\vt\|_{H^2_{tan}(\Gamma)}\|\ut\|_{H^2_{tan}(\Gamma)}.\label{eq:geomA2}
 \end{align}
 
\end{lemma}

\begin{proof}
\textbf{Estimate (\ref{eq:geomF})} is proven by applying the change of variable. Let 
{the surface measures of $\Gamma$ and $\Gamma_h$ be given by $d\sigma$ and $d \sigma_h$, respectively.  Following \cite{Demlow2007}, for any $x \in \Gamma_h$ there exists an $\nu_h$ such that $d\sigma_h \circ p(x)= \nu_h(x) d \sigma_h(x)$ and the following estimate holds: } 
\begin{equation}\label{eq:nu}
\|1-\nu_h\|_{L^\infty(K)}\leq ch^2. 
\end{equation}
We get 
\begin{align*}
   | (\ft,\vt_h^l)_\Gamma-(\tilde \ft,\vt_h)_{\Gamma_h}|\leq |((\nu_h-1)\tilde \ft,\vt_h)_{\Gamma_h}| \leq c h^2 \|\ft\|_{L^2(\Gamma_h)} \|\vt\|_{L^2(\Gamma_h)}.
\end{align*}
\textbf{Estimate (\ref{eq:geomA})} and \textbf{Estimate (\ref{eq:geomA2})}  are proven by applying estimate (5.31) and estimate (5.33) of Lemma 5.5 in \cite{HansboLarson2020} with $k_g=1$, respectively.
\end{proof}
 {
 \subsection{Bound for the cross-element conormal jump}\label{sec:appsmall}
 \begin{lemma}[Bound for the cross-element conormal jump]\label{lemma:jumpSmall}
 Let $\ut$ be the solution of (\ref{eq:ana}) and $ \tilde \ut$ be its extension to $U$ defined by (\ref{eq:extension}). Then the following estimate holds
 \begin{equation}
 \Big |\int_{E^l}[v_n^l]\nt^+_{E^l} \cdot \nabla_\Gamma \ut \nt^+_{E^l}\, d\sigma \Big |\\ \leq h \|\vt_h\|_{h}\|\ut\|_{H^2(K^l)}\quad \forall \vt_h\in \Vt_h
 \end{equation}
 where $[v_n]:=v^+_n + v^-_n$  with $v^{+}_n=\vt_h|_{K^+} \cdot \nt^{+}_{E}$ and $v^{-}_n=\vt_h|_{K^-} \cdot \nt^{-}_{E}$ defines a jump based on $\nt_E^+$ and $\nt^-_E$.
 \end{lemma}
 \begin{proof}
Let ${\bar v_n^{l,+}=\frac{1}{|E^l|}\int_{E^l}} v_n^{l,+}\, d\sigma$ and ${\bar v_n^{-,l}=\frac{1}{|E^l|}\int_{E^l}} v_n^{l,-}\, d\sigma$ denote the averages of $v^{l,+}_n$ and  $v^{l,-}_n$  over an edge $E^l$, respectively. We introduce the mean value over the jump and get
 \begin{equation}\label{eq:jump3}
\begin{aligned}
\int_{E^l}[v_n^l]&\nt^+_{E^l}\cdot \nabla_\Gamma \ut \nt^+_{E^l} \,d\sigma \\
&=\int_{E^l} \Big ([v_n^l ]-\overline{[ v_n^l]}\Big)\nt^+_{E^l} \cdot \nabla_\Gamma \ut \nt^+_{E^l} \,d\sigma  + \int_{E^l}  \overline{[v_n^l]}\nt^+_{E^l} \cdot \nabla_\Gamma \ut\nt^+_{E^l} \,d\sigma \\
&=\int_{E^l} \Big ([v_n^l - \bar v_n^l]\Big){(\nt^+_{E^l})^T}   \nabla_\Gamma \ut \nt^+_{E^l} \,d\sigma + \int_{E^l} [\bar v_n^l]\nt^+_{E^l} \cdot \nabla_\Gamma \ut \nt^+_{E^l}\, d\sigma\\
&=I^1_{E^l}+I^2_{E^l}. 
\end{aligned}
\end{equation}
We treat  the resulting terms in (\ref{eq:jump3}) separately.

\textbf{Term $I^1_{E^l}$} is bounded by  using the fact that \begin{equation}
\int_E  [v_n - \bar v_n] =0, 
\end{equation} and introducing a constant in form of the mean value $\bar \rho$ with $\rho:={(  \nt^+_{E^l})^T} \nabla_\Gamma \ut \nt^+_{E^l}$. To improve readability, we use the subscript e to indicate the extended function. 
 \begin{align}
\int_{E^l}& \Big ([v_n^l - \bar v_n^l]\Big) {(\nt^+_{E^l})^T}   \nabla_\Gamma \ut \nt^+_{E^l} \,d\sigma \\
\leq& \Big |\int_{E}[v_n-\bar v_n] \cdot \Big ( \rho^e  -   (\bar\rho)^e \Big) \nu_hd\sigma_h \Big |\\ 
\leq& \Big (\int_{E} |[v_n - \bar v_n]|^2 \,d\sigma_h \Big)^\frac{1}{2} \Big (\int_{E^l} |\rho-\bar \rho|^2\, d\sigma \Big)^\frac{1}{2},
\end{align}
where we use the Cauchy-Schwarz inequality in the last estimate. 
We apply the discrete Poincare inequality and the trace inequality and get
\begin{equation}\label{eq:jump6}
\begin{aligned}
 \Big (\int_{E} |[v_n - \bar v_n]|^2 \,\nu_hd\sigma_h \Big)^\frac{1}{2}  &\leq ch^\frac{1}{2} \Big (\|\vt_h^+\|_{H^1(K^+)}+  \|\vt_h^-\|_{H^1(K^-)} \Big),\\
 \Big (\int_{E} |\rho-\bar\rho |^2\, d\sigma \Big)^\frac{1}{2} &\leq c h^\frac{1}{2} \|\ut\|_{H^2(K^l)}.
\end{aligned}
\end{equation} 

\textbf{Term $I^2_{E^l}$}  is bounded by applying the Cauchy-Schwarz inequality  
\begin{equation}
\begin{aligned}
\int_{E^l}& [\bar v_n^l]\nt^+_{E^l} \cdot \nabla_\Gamma \ut \nt^+_{E^l}\, d\sigma\\
 &\leq \Big(\int_{E^l} |[\bar v_n^l]\nt^+_{E^l}|^2 \Big)^\frac{1}{2} \Big(\int_{E^l} |\nabla_\Gamma \ut \nt^+_{E^l}|^2 \, d\sigma \Big)^\frac{1}{2} \\
&\leq c {h} \Big( \|\vt_h\|_{h(K^-)}+\|\vt_h\|_{h(K^+)}\Big )\|\ut\|_{ H^2(K^l)}.
\end{aligned}
\end{equation}
The bound for $I^2_{E^l}$ results by applying the  trace inequality and noting that 
 \begin{equation}
\int_{E^l} [\bar v^l_n] \,d\sigma = \int_{E} [\bar v_n ]\nu_h \,d\sigma_h = (1-\nu_h) \int_E [\bar v_n]\,d\sigma_h +\underbrace{ \int_E [\bar v_n]\,d\sigma_h}_{=0} \leq  c h^2 \int_E [\bar v_n]\,d\sigma_h,
 \end{equation}
 where we use the estimate (\ref{eq:nu}) in the last bound.
  Thus, we get 
\begin{equation}\label{eq:jump5}
 I^2_{E^l}\leq c h\Big( \|\vt_h\|_{h(K^-)}+\|\vt_h\|_{h(K^+)}\Big )\|\ut\|_{ H^2(K^l)}.
\end{equation}
Combining (\ref{eq:jump3}) with (\ref{eq:jump5}) and  (\ref{eq:jump6}) yields the final estimate.
 \end{proof}

 \subsection{Proof to Lemma \ref{lemma:Klaus2}}\label{sec:appbig}
   \begin{proof}
 We note that $\Delta_{\Gamma_h}({\Pi^{tan}_h}\ut)=0$. By including the finite-element interpolation $ \Pi_h^{tan} \ut $ and applying Green’s first identity (\ref{eq:green}), we obtain
 \begin{equation}
 \begin{aligned}
       a_h(\tilde \ut,\tilde \zt- {\Pi^{tan}_h} \zt)
  =&a_h({\Pi^{tan}_h} \ut,\tilde \zt-{\Pi^{tan}_h} \zt )+a_h(\tilde \ut-{\Pi^{tan}_h} \ut,\tilde \zt-{\Pi^{tan}_h}  \zt)\\
   =&\sum_K
  \Big ({\Pi^{tan}_h} \ut -(\underbrace{\Delta_{\Gamma_h}( {\Pi^{tan}_h} \ut)}_{=0}, \tilde \zt -{\Pi^{tan}_h} \zt \Big )_K\\
   +& \sum_{E}\underbrace{ \Big ( \nabla^+_{\Gamma_h} ({\Pi^{tan}_h} \ut) \nt^+_E+\nabla^-_{\Gamma_h} ({\Pi^{tan}_h} \ut)\nt^-_E, \tilde \zt-{\Pi^{tan}_h} \zt\Big )_E}_{I^E_1}\\
  +&\sum_K
  \Big (\tilde \ut -{\Pi^{tan}_h} \ut -(\Delta_{\Gamma_h}( \tilde \ut -{\Pi^{tan}_h} \ut), \tilde \zt -{\Pi^{tan}_h} \zt \Big )_K\\
  +& \sum_{E}\underbrace{ \Big ( \nabla^+_{{\Gamma_h}} ( \tilde \ut-{\Pi^{tan}_h}) \nt^+_E+\nabla^-_{\Gamma_h} (\tilde \ut-{\Pi^{tan}_h} \ut)\nt^-_E, \tilde \zt- {\Pi^{tan}_h} \zt\Big )_E}_{I^E_2}.
 \end{aligned}
 \end{equation}
Applying the Cauchy–Schwarz inequality together with $L^2$-interpolation estimate (\ref{eq:L2est}) and (\ref{eq:L2estPh}), we obtain 
 \begin{equation}\label{eq:PrimDualEq}
 \begin{aligned}
       a_h(\tilde \ut,\tilde \zt-{\Pi^{tan}_h} \zt)
  \leq&\sum_K \Big( \int_K |\Delta_{\Gamma_h} \tilde \ut|^2  d\sigma_h\Big)^\frac{1}{2}\Big(\int_K | \Pt_h(\tilde \zt- {\Pi^{tan}_h}\zt) |^2  d\sigma_h\Big)^\frac{1}{2}\\
  &+\sum_K
\Big ( \int_K |\tilde \ut|^2  d\sigma_h \Big )^\frac{1}{2} \Big( \int_K |\Pt(\tilde \zt-{\Pi^{tan}_h}) \zt|^2  d\sigma_h \Big)^\frac{1}{2} +I^E_1, +I^E_2\\
  \leq&  \sum_{K}
 c h^2 \|\ut\|_{H^2(K)} \|\zt\|_{H^2(K)} +I^E_1+I^E_2.
 \end{aligned}
 \end{equation}
 
\textbf{Term $I^E_1$:} We
apply the Cauchy-Schwarz inequality and get
\begin{align}\label{eq:I1}
I^E_1= \underbrace{\Big (\int_E | \nabla^+_{\Gamma_h} ({\Pi^{tan}_h} \ut) \nt^+_E+\nabla^-_{\Gamma_h} ({\Pi^{tan}_h} \ut)\nt^-_E|^2\, d\sigma_h\Big)^\frac{1}{2}}_{I^E_{11}}\underbrace{\Big( \int_E| \Pt_h(\tilde \zt- {\Pi^{tan}_h} \zt)|^ 2\, d\sigma_h\Big )^\frac{1}{2}}_{\leq c h^\frac{3}{2} \|\tilde \zt\|_{H^2(K)}, (\ref{eq:L2estPurPh})}.
\end{align}
To bound $I^E_{11}$ we introduce the  covariant derivative on $\Gamma$, and use the fact that $\nt^+_{E^l}=-\nt^-_{E^l}$. To improve readability, we indicate the extended functions by an $e$ in the following estimate.
\begin{align*}
I^E_{11}
\leq& 
|\nabla^+_{\Gamma_h} {\Pi^{tan}_h} \ut- (\nabla_\Gamma ({\Pi^{tan}_h} \ut)^l )^e|_{L^2(E)}|\nt_E^+|_{L^2(E)}\\
+& |\nabla^-_{\Gamma_h} {\Pi^{tan}_h} \ut-({\nabla_\Gamma ({\Pi^{tan}_h} \ut)^l})^e|_{L^2(E)}|\nt_E^-|_{L^2(E)}\\
+&|({\nabla_\Gamma ({\Pi^{tan}_h} \ut)^l})^e|_{L^2(E)}\Big(|\nt_E^+ -\nt^+_{E^l}|_{L^\infty(E)}+|\nt_E^--\nt^-_{E^l}|_{L^\infty(E)}\Big).
\end{align*}
We apply the trace inequality, the estimate to bound the operator differences given in Lemma \ref{lemma:OpDiff}, the $H^1$-stability of the interpolation and the geometric bound (\ref{eq:nel}) and obtain 
\begin{align}\label{eq:boundI11}
I^E_{11}\leq& { c h^\frac{1}{2}\|\tilde \ut \|_{H^2(K)}}.
\end{align}
Employing (\ref{eq:boundI11}) to estimate (\ref{eq:I1}), we arrive at the following bound for $I_1^E$:
\begin{align}\label{eq:estI1E}
I^E_{1}\leq c h^\frac{1}{2}\|\tilde \ut\|_{H^2(K)} h^\frac{3}{2}\|\tilde \zt\|_{H^2(K)}.
\end{align}
\textbf{Term $I^E_2$:}
We use the Cauchy-Schwarz inequality and get 
\begin{align*}
  I^E_2\leq I_{21}^\frac{1}{2} I_{22}^\frac{1}{2},\\
  I_{21}&=\int_E| \nabla_{\Gamma_h}^+ (\tilde \ut-{\Pi^{tan}_h} \ut) \nt_E^+ + \nabla_{\Gamma_h}^-(\tilde \ut-{\Pi^{tan}_h} \ut) \nt_E^-|^2  d\sigma_h,  \\
  I_{22}&=\int_E| \Pt_h(\tilde \zt -{\Pi^{tan}_h} \zt) |^2  d\sigma_h.
\end{align*}
Note that $E$ is the edge of the flat triangle $K$. Based on the trace inequality, the interpolation estimate (\ref{eq:H1est}), and the fact that $\nabla_{\Gamma_h}^2( \Pi_K^{tan} \vt)=0$, we bound $I_{21}$ as
\begin{equation}            
\begin{aligned}
 I_{21}
   &  \leq  h\|\tilde \ut\|^2_{H^2(K^-)}+ h\|\tilde \ut\|^2_{H^2(K^+)}.
     \end{aligned}
\end{equation}
A bound for  $I_{22}$ is obtained  by applying the interpolation estimate (\ref{eq:L2estPurPh}). 
This gives us
\begin{align}\label{eq:I2E}
    I^E_2\leq c h^\frac{1}{2}\|\tilde \ut\|^2_{H^2(K)}  h^\frac{3}{2}\|\tilde \zt\|^2_{H^2(K)}.
\end{align}

 Finally, we obtain the desired estimate by combining (\ref{eq:PrimDualEq}) with (\ref{eq:estI1E}) and (\ref{eq:I2E}), and applying the norm equivalence (\ref{eq:norm2}) and using the regularity estimate (\ref{eq:regu}).
\begin{align*}
     a_h(\tilde \ut,\tilde \zt -\Pt{\Pi^{tan}_h}\zt)\leq c h^2 \|\ft\|_{L^2(\Gamma)}\|\zt\|_{H^2(\Gamma)}.
\end{align*}
 \end{proof}

\end{document}